\documentclass[a4paper,reqno,10pt,oneside]{amsart}
\usepackage{amsmath}
\usepackage[left=2.61cm,right=2.61cm,top=2.72cm,bottom=2.72cm]{geometry}
\usepackage[colorlinks=true,urlcolor=blue,
citecolor=red,linkcolor=blue,linktocpage,pdfpagelabels,
bookmarksnumbered,bookmarksopen]{hyperref}
\usepackage{graphicx}
\usepackage{mathrsfs}
\usepackage{amssymb,amsmath, amsfonts, amsthm}
\usepackage{latexsym}

\allowdisplaybreaks[1]

\numberwithin{equation}{section}

\renewcommand{\(}{\left(}
\renewcommand{\)}{\right)}
\renewcommand{\[}{\left[}
\renewcommand{\]}{\right]}

\newtheorem{theorem}{Theorem}[section]
\newtheorem{proposition}[theorem]{Proposition}

\newtheorem{lemma}[theorem]{Lemma}
\theoremstyle{definition}
\newtheorem{remark}[theorem]{Remark}
\theoremstyle{definition}
\newtheorem{definition}[theorem]{Definition}
\theoremstyle{definition}

\renewcommand{\le}{\leqslant}
\renewcommand{\ge}{\geqslant}

\renewcommand{\l }{\langle}
\renewcommand{\r }{\rangle}

\renewcommand{\S}{\Sigma}

\newcommand{\T}{{\mathcal T}}

\renewcommand{\S}{\mathbb{S}}

\newcommand{\beq}{\begin{equation}}
\newcommand{\eeq}{\end{equation}}
\newcommand{\beqs}{\begin{equation*}}
\newcommand{\eeqs}{\end{equation*}}
\newcommand{\beqn}{\begin{eqnarray}}
\newcommand{\eeqn}{\end{eqnarray}}
\newcommand{\beqns}{\begin{eqnarray*}}
\newcommand{\eeqns}{\end{eqnarray*}}
\newcommand{\bdoc}{\begin{document}}
\newcommand{\edoc}{\end{document}}
\newcommand{\be}{\begin{enumerate}}
\newcommand{\ee}{\end{enumerate}}
\newcommand{\bdescr}{\begin{description}}
\newcommand{\edescr}{\end{description}}
\newcommand{\ba}{\begin{array}}
\newcommand{\ea}{\end{array}}
\newcommand{\intR}{\int_{\mathbb R^N}}
\newcommand{\R}{\mathbb R}
\newcommand{\B}{\mathbb B}
\newcommand{\C}{\mathbb C}
\renewcommand{\H}{\mathbb H}
\renewcommand{\L}{\mathbb L}
\newcommand{\parallelsum}{\mathbin{\!/\mkern-5mu/\!}}
\newcommand{\e}{\epsilon}
 \renewcommand{\(}{\left(}
\renewcommand{\)}{\right)}
\renewcommand{\[}{\left[}
\renewcommand{\]}{\right]}
\newenvironment{Proof}{\noindent{\bf Proof}}{\hfill$\Box$\\[2mm]}

\begin{document}
\title[Radial graphs]{Spacelike radial graphs of prescribed mean curvature\\ in the Lorentz-Minkowski space}

\author{D\MakeLowercase{enis} Bonheure, A\MakeLowercase{lessandro} Iacopetti}

\thanks{\emph{Acknowledgements.} Research partially supported by the project ERC Advanced Grant 2013 n.~339958 Complex Patterns for Strongly Interacting Dynamical Systems COMPAT, by Gruppo Nazionale per l'Analisi Matematica, la Pro\-ba\-bi\-li\-t\`a e le loro Applicazioni (GNAMPA) of the Istituto Nazionale di Alta Matematica (INdAM) by FNRS (PDR T.1110.14F and MIS F.4508.14) and by ARC AUWB-2012-12/17-ULB1- IAPAS}
\address[Denis Bonheure]{D\'epartement de math\'ematique, Universit\'e Libre de Bruxelles, Campus de la Plaine - CP214 boulevard du Triomphe, 1050 Bruxelles, Belgium}
\email{denis.bonheure@ulb.ac.be}
\address[Alessandro Iacopetti]{D\'epartement de math\'ematique, Universit\'e Libre de Bruxelles, Campus de la Plaine - CP214 boulevard du Triomphe, 1050 Bruxelles, Belgium}
\email{alessandro.iacopetti@ulb.ac.be}
\begin{abstract}

In this paper we investigate the existence and uniqueness of spacelike radial graphs of prescribed mean curvature in the Lorentz-Minkowski space $\L^{n+1}$, for $n\geq 2$, spanning a given boundary datum lying on the hyperbolic space $\H^n$.
\end{abstract}

\maketitle

\section{Introduction}

A radial graph is a hypersurface $\Sigma$ such that each ray emanating from the origin intersects $\Sigma$ once at most. In the euclidean context the problem of finding radial graphs of prescribed mean curvature has been extensively studied over the years. The first paper on the subject is due to Rad\`o. In his paper \cite{Rado}, he proved that for any given Jordan curve  $\Gamma \subset \R^3$, with one-one radial projection onto a convex subset of the unit sphere $\S^2$, there exists a minimal graph spanning $\Gamma$. Later, Tausch proved, in \cite{Tausch}, that area-minimizing disk-type hypersurfaces spanning a boundary datum $\Gamma$ which can be expressed as a radial graph over $\partial\Omega$, where $\Omega \subset \S^n$ is a convex subset, have a local representation as a radial graph. The case of variable mean curvature was investigated by Serrin \cite{Serrin2}, and a recent result of radial representation for $H$-surfaces in cones has been given in \cite{CaldIac}. Treibergs and Wei studied in \cite{TW} the case of closed hypersurfaces, i.e. compact hypersurfaces without boundary. Lopez in \cite{Lo03} and de Lira in \cite{HerbDeLira} studied the case of radial graphs of constant mean curvature.\\ 

The Lorentz-Minkowski space, denoted by $\L^{n+1}$, is defined as the vector space $\R^{n+1}$ equipped with the symmetric bilinear form
$$\langle x,y \rangle := x_1y_1+\ldots+x_ny_n-x_{n+1}y_{n+1},$$
where $x=(x_1,\ldots,x_{n+1}), y=(y_1,\ldots,y_{n+1}) \in \R^{n+1}$.
The  bilinear form $\langle \cdot, \cdot \rangle$ is a non-degenerate bilinear form of index one (see \cite[Sect. A]{Spivak}), where the index of a bilinear form on a real vector space is defined as the largest dimension of a negative definite subspace. The modulus of $v \in \L^{n+1}$ is defined as $|v|:=\sqrt{|\l v,v\r|}$.\\

The interest on finding spacelike hypersurfaces of prescribed mean curvature in the Lorentz-Minkowski space comes from the theory of relativity, in which maximal and constant mean curvature spacelike hypersurfaces play an important role (see \cite{{BartnikSimon}}), where spacelike means that the restriction of the Lorentz metric to the tangent plane, at every point, is positive definite. In the literature, several result are available for spacelike vertical graphs, i.e. hypersurfaces which are expressed as a cartesian graph. Entire maximal spacelike hypersurfaces were studied by Cheng and Yau in \cite{ChengYau} and later Treibergs  in \cite{T82} tackled the general case of entire spacelike hypersurfaces of constant mean curvature. The Dirichlet problem for spacelike vertical graphs in $\L^{n+1}$ was solved by Bartik and Simon in \cite{BartnikSimon}, and Gerhardt in \cite{Gerhardt} extended those results to the case of vertical graphs contained in a Lorentzian manifolds which can be expressed as a product of a Riemannian manifold times an interval. Bayard studied in \cite{Bayard2} the more general problem of prescribed scalar curvature. On the contrary, for radial graphs, to our knowledge, the only available result concerns entire spacelike hypersurfaces with prescribed scalar curvature which are asymptotic to the light-cone (see \cite{Bayard1}).\\

The geometry of Lorentz-Minkowski spaces plays an in important role in the setting of the problem. A first relevant fact is that there cannot exist spacelike closed hypersurfaces (see Proposition \ref{nonexistclosedspacelike}, or \cite{Lopez2014} for the case of surfaces in $\L^3$). Therefore $\S^n$-type surfaces are ruled out, and the model hypersurface in $\L^{n+1}$ for describing spacelike radial graphs is the hyperbolic space $\H^n$ (see Definition \ref{defhyperbolicspace}). Another important feature of Lorentz-Minkowski spaces is that, given a domain, there exist spacelike hypersurfaces of arbitrarily large (in modulus) mean curvature (see \cite{Lopez2013}), while in the euclidean context this is not true in general. This fact will be crucial in our paper to construct barriers.\\

We state now the problem. Let $\Omega$ be a smooth bounded domain of $\H^n$. For $u:\overline\Omega \to \R$, we define the associated radial graph over $\Omega$ as the set $$\Sigma(u):=\{p=e^{u(q)}q \in \L^{n+1}; \ q \in \overline\Omega\}.$$ 
Let  $\mathcal{C}_{\overline\Omega}$ be the cone spanned by $\overline\Omega$ (minus the origin), i.e.
$\mathcal{C}_{\overline\Omega}:=\{p=\rho q \in \L^{n+1}; \ q\in \overline\Omega,\ \rho >0\},$
and let $H: \mathcal{C}_{\overline\Omega} \to \R$. 

\begin{definition}
A $H$-bump (over $\Omega$) is a radial graph $\Sigma$ whose boundary coincides with $\partial \Omega$, and such that the mean curvature of $\Sigma$ at every (interior) point equals $H$. 
\end{definition}

The Dirichlet problem for spacelike $H$-bumps is given by

 \begin{equation}\label{eqHbump}
\begin{cases}
 \displaystyle \sum_{i,j=1}^n \left((1-|\nabla u|^2)\delta_{ij} + u_i u_j)\right)u_{ij}= n (1-|\nabla u|^2) - n  (1-|\nabla u|^2)^{3/2}e^uH(e^uq) & \hbox{in} \ \Omega,\\[12pt]
|\nabla u|<1  &  \hbox{in} \ \Omega,\\[6pt]
u=0 & \hbox{on} \ \partial\Omega,
\end{cases}
\end{equation}
 where $u_i$, $u_{ij}$ are the covariant derivatives of $u$, $\nabla u$ is the gradient with respect to the Levi-Civita connection of $(\H^n,g)$ (see Sect. 3), and $g=dx_1 \otimes dx_1+ \ldots +dx_n \otimes dx_n - dx_{n+1} \otimes dx_{n+1}$ is the induced Riemannian metric on $\H^n$ (see Sect. 2).

\begin{definition}\label{defconicalcap}
Let $0<r_1 \leq 1 \leq r_2$, with $r_1\neq r_2$. The hyperbolic conical cap of radii $r_1,r_2$ spanned by $\overline\Omega$ is the set  $\mathcal{C}_{\overline{\Omega}}(r_1,r_2):=\{p=\rho q \in \L^{n+1}; \ q\in \overline{\Omega},\ r_1 \leq \rho \leq r_2\}.$
\end{definition}

 The main result of our paper is the following existence theorem.

 \begin{theorem}\label{mainteoprob}
Let $\alpha \in (0,1)$, $0<r_1 \leq 1 \leq r_2$, with $r_1\neq r_2$. Assume $\Omega$ is a bounded domain of $\H^n$ of class $C^{3,\alpha}$ that satisfies a uniform exterior geodesic ball condition. If $H \in C^{1,\alpha}({C_{\overline{\Omega}}(r_1,r_2)})$ is positive and satisfies
\begin{itemize}
\item[i)] $H(r_1 q) > r_1^{-1} \ \hbox{and} \ H(r_2 q) < r_2^{-1} \ \hbox{for any} \ q \in \overline{\Omega}$;
\item[ii)] $\frac{\partial}{\partial \lambda} \left(\lambda H(\lambda q)\right) \leq 0$, for all $q \in {\overline\Omega}$, $\lambda \in [r_1,r_2]$;
\end{itemize}
then there exists a unique solution of Problem \ref{eqHbump} whose associated radial graph is contained in $\mathcal{C}_{\overline\Omega}(r_1,r_2)$. 
\end{theorem}

Let $\Omega$, $r_1,r_2$ be in the statement of Theorem \ref{mainteoprob}. Let $m\geq1$, let $\omega:\overline\Omega \to \R^+$ be a smooth positive function such that $r_1^{m-1}<\omega<r_2^{m-1}$ and let $H_{m,\omega}:\mathcal{C}_{\overline\Omega}(r_1,r_2) \to \R^+$, defined by 
\begin{equation}\label{modelfunctionH}
H_{m,\omega}(x):=\frac{\omega\left(\frac{x}{|x|}\right)}{\ \ |x|^m}.
\end{equation}
One easily verifies that $H_{m,\omega}$ satisfies the hypotheses ${\rm i)}$ and ${\rm ii)}$ of Theorem \ref{mainteoprob}.
In particular, this shows the existence of spacelike radial graphs of prescribed mean curvature even for non homogenous functions $H$, a case which is not contemplated for instance in \cite{Bayard1}, where the $k$-th scalar curvature is prescribed just on $\H^n$.\\

We remark that, the equation \eqref{eqHbump} can be put in divergence form, namely
\begin{equation}\label{eqdivform}
\begin{cases}
\displaystyle -\mathrm{div}_{\H^n}\left(\frac{\nabla u}{\sqrt{1-|\nabla u|^{2}}}\right)+\frac{n}{\sqrt{1-|\nabla u|^{2}}}=n e^{u}H(e^{u}q)&\text{in }\Omega,\\
|\nabla u|<1&\text{in }\Omega,\\
u=0&\text{on }\partial\Omega,
\end{cases}
\end{equation}
where $\mathrm{div}_{\H^n}$ denotes the divergence operator for $(\H^n,g)$. The principal part of this operator appears in the Born-Infeld theory of electromagnetism \cite{BI-proc}, which is a particular example of what is usually known as a nonlinear electrodynamics. We therefore stress that Theorem \ref{mainteoprob} provides existence and uniqueness of solutions for some specific Born-Infeld equations in which appear non trivial nonlinearities involving both the gradient and the function (see also \cite{BDP, BI-ARMA} and the references therein).\\

The proof of Theorem \ref{mainteoprob} relies on the combinations of several tools. For the existence, we apply a variant of the classical Leray-Schauder fixed point theorem due to Potter, see \cite{Potter}. To this aim, we make use of suitable comparison theorems and we prove fine a priori estimates for the solutions and their gradient. Regarding uniqueness, we take benefit of the Hopf maximum principle as in the version stated by Pucci and Serrin in \cite{PuSe04}. 

We point out that the uniform exterior geodesic ball condition allows to construct barriers for the gradient of the solutions at the boundary. Such construction strongly depends on the shape of the mean curvature operator for spacelike hypersurfaces in the Lorentz-Minkowski space, and we remark that Theorem \ref{mainteoprob} grants existence of spacelike radial graphs over arbitrary large and even non convex domains of $\H^n$. We note that it is not possible to mimic this construction in the euclidean framework, and in fact the problem of finding radial graphs over proper (possibly non convex) domains of $\S^n$ which are not contained in a hemisphere is still open.\\

Concerning global a priori estimates for the gradient, which is the key step in the proof, we derive a quite complex technical result, see Proposition \ref{fineraprioriestgrad}, which is inspired from the paper \cite{Gerhardt} and is based on the introduction of an ad hoc differential operator, Stampacchia's truncation method and fine estimates of the $L^p$-norm of the quantity $\nu(u)=\frac{1}{\sqrt{1-|\nabla u|^2}}$.\\

In this paper we also introduce a new definition of admissible couple $(\Omega,H)$ and triple $(\Omega,H,\theta)$, see Definition \ref{classadmissible}, where $\theta \in (0,1)$. This notion of admissibility is very general and works even for non smooth domains and just for continuous functions $H$. However, given a couple $(\Omega,H)$, it is not easy in general to verify whether it is admissible or not. In Sect. 4 we provide trivial examples of admissible couples and in Proposition \ref{propadmissextbgball} we exhibit a class of functions $H$ such that $(\Omega,H)$ is admissible whenever $\Omega$ satisfies a uniform exterior geodesic condition. Using the notion of admissible couple, we can extend Theorem \ref{mainteoprob} to a wider class of domains and mean curvature functions.

 \begin{theorem}\label{mainteoprob2}
Let $\alpha \in (0,1)$, $0<r_1 \leq 1 \leq r_2$, with $r_1\neq r_2$. Assume that $\Omega$ is a bounded domain of $\H^n$ 
of class $C^{3,\alpha}$ and $H \in C^{1,\alpha}({C_{\overline{\Omega}}(r_1,r_2)})$ satisfies the conditions {\rm i)} and {\rm  ii)} of Theorem \ref{mainteoprob}. Assume that $(\Omega,H)$ is admissible. Then there exists a unique solution of Problem \ref{eqHbump} whose associated radial graph is contained in $\mathcal{C}_{\overline\Omega}(r_1,r_2)$. 
\end{theorem}

A further existence result for Problem \ref{eqHbump}, under more restrictive assumptions, is as follows. 

 \begin{theorem}\label{mainteo}
Let $\alpha \in (0,1)$ and $\Omega$ be a bounded domain of $\H^n$ of class $C^{3,\alpha}$. Assume  $\theta\in(0,1)$, $0<r_1 \leq 1 \leq r_2$, with $r_1\neq r_2$, and   $H \in C^{1,\alpha}(\mathcal{C}_{\overline{\Omega}}(r_1,r_2))$ satisfies
\begin{itemize}
\item[a)] $H(r_1 q) > r_1^{-1} \ \hbox{and} \ H(r_2 q) < r_2^{-1} \ \hbox{for any} \ q \in \overline{\Omega}$;
\item[b)] $\frac{\partial}{\partial \lambda} \left(\lambda H(\lambda q) \right)< - \frac{1}{r_1(\theta-\theta^2/4)^{1/2}}$, for all $q \in \overline\Omega$, $\lambda \in [r_1,r_2]$;
\item[c)] $\left\|\nabla_0^T H(x)\right\|_{n+1} < \frac{1-\theta}{n^{3/2}r_2^2}$, for all  $x\in \mathcal{C}_{\overline{\Omega}}(r_1,r_2)$, where $\nabla_0^T H$ is the euclidean tangential component of $\nabla_0 H(x)$ on $T_{x/|x|} \H^n$ (see Definition \ref{tangentialcomponent}), $\nabla_0H$ is the gradient of $H$ with respect to the euclidean flat metric, $\|\cdot\|_{n+1}$ is the euclidean norm in $\R^{n+1}$.
\end{itemize}
Assume at last that $(\Omega, H,\theta)$ is admissible according to Definition \ref{classadmissible} and Definition \ref{classadmissibleext}. Then there exists a unique spacelike $H$-bump contained in $\mathcal{C}_{\overline{\Omega}}(r_1,r_2)$. 
\end{theorem}

We mention this result because the proof quite differs from that of Theorem \ref{mainteoprob2} and better enlightens the differences and difficulties with respect to the euclidean case. The proof is this time based on the classical Leray-Schauder theorem, see for instance \cite[Theorem 11.3]{GT}. The first step is to solve a suitable regularized equation associated to \eqref{eqHbump}, see \eqref{eqpseudoHbumps} and Theorem \ref{teopseudosol}. The idea of solving such a regularized equation is taken from \cite{T82} where the author construct barriers for the gradient at the boundary. The way back to the original Dirichlet problem then uses a gradient maximum principle \cite[Proposition 6]{TW}. In contrast with \cite{T82}, we deal here with equations which do not satisfy, in general, a gradient maximum principle \cite[Theorem 15.1]{GT}. In fact, in our case, when passing to local coordinates, we see that the regularized operator associated to \eqref{eqHbump} does not satisfy, in general \cite[condition (15.11)]{GT}, and the principal part depends both on the gradient and on the domain variables. We refer to Lemma \ref{lemmaunifelliptic} below for more details. In order to overcome this difficulty, and eventually deduce a global a priori $C^1$ estimate, we perform the regularization in a proper way. We then use the admissibility condition to control the gradient at the boundary, whereas we use two different strategies, see Lemma \ref{lemapriorigradestimates}, for the interior estimate. The first one which is based on the properties of harmonic functions, works only in dimension two. The other proof works in any dimension and is based on the global gradient bound given by  \cite[Theorem 15.2]{GT}. 

At last, in the spirit of \cite{TW}, we prove a new kind of interior gradient estimate, see Proposition \ref{stimagradint}, so that, under the hypotheses of Theorem \ref{mainteo}, the solution of the regularized problem is a solution of \eqref{eqHbump}. It is important to note that, in contrast to \cite{TW}, since $\H^n$ has negative Ricci curvature and since we deal with hypersurfaces with boundary, the mere gradient estimate of Proposition \ref{stimagradint} is not sufficient for getting a global a priori $C^1$-estimate. We refer to Remark \ref{remsustimagrad} for more details.

\medbreak

When $\Omega$ satisfies a uniform exterior geodesic conditions, thanks to Proposition \ref{propadmissextbgball}, Remark \ref{rmk:thetadipendenceH} and Remark \ref{remextensadmissible}, it is possible to show that the functions given by \eqref{modelfunctionH} satisfy the hypotheses of Theorem \ref{mainteo} for suitable choices of $r_1,r_2,m$, for $\omega$ close to one (in the $C^1$-topology), and for some $\theta_* \in (0,1)$. 

\medbreak

As a future goal, it would be natural  to investigate if it is possible to remove the monotonicity assumption on $H$ and to extend Theorem \ref{mainteoprob} also to sign-changing mean curvature functions. 

\medbreak

The outline of the paper is the following. In Sect. 2, we fix the notations and we collect some known facts which are useful in the remaining of the paper. In Sect. 3, we derive the equation for spacelike $H$-bumps and in Sect. 4 we prove Proposition \ref{propadmissextbgball} and some a priori estimates. Sect. 5 is dedicated to the proof of existence and uniqueness of solutions for the regularized Dirichlet problem associated to Problem \ref{eqHbump}. In Sect. 6, we work out an interior gradient estimate, namely Proposition \ref{stimagradint}, and in Sect. 7 we prove Theorem \ref{mainteo}. In Sect. 8, we prove a global a priori estimate for the gradient. We finally prove Theorem \ref{mainteoprob} and Theorem \ref{mainteoprob2} in Sect. 9.

\section{Notation and preliminary results}
Let $n\geq 2$, we denote by $\L^{n+1}$ the (n+1)-dimensional Lorentz-Minkowski space, which is $\R^{n+1}$ equipped with the symmetric bilinear form 
$$\langle x,y \rangle := x_1y_1+\ldots+x_ny_n-x_{n+1}y_{n+1}.$$
We classify the vectors of $\L^{n+1}$  in three types.

\begin{definition}\label{deftipology}
A vector $v \in \L^{n+1}$ is said to be
\begin{itemize}
\item spacelike if $\l v,v\r>0$ or $v=0$;
\item timelike if $\l v,v\r<0$;
\item lightlike if $\l v,v\r=0$ and $v \neq 0$.
\end{itemize}
\end{definition}
The modulus of $v \in \L^{n+1}$ is defined as $|v|:=\sqrt{|\l v,v\r|}$. We also denote by $(x,y)_{n+1}= x_1y_1+ \ldots + x_{n+1}y_{n+1}$ the euclidean scalar product, and by $\|x\|_{n+1}=\sqrt{x_1^2+\ldots+x_{n+1}^2}$ the euclidean norm in $\R^{n+1}$. 
Given a vector subspace $V$ of $\L^{n+1}$, we consider the induced metric $\l\cdot,\cdot\r_V$ defined in the natural way
$$ \l v,w\r_V:=\l v, w \r, \ v,w \in V.$$
According to Definition \ref{deftipology} we classify the subspaces of $\L^{n+1}$ as follows.
\begin{definition}
 A vector subspace $V$ of $\L^{n+1}$ is said to be:
\begin{itemize}
\item spacelike if the induced metric is positive definite;
\item timelike if the induced metric has index 1;
\item lightlike if the induced metric is degenerate.
\end{itemize}
\end{definition}

In this paper, we deal only with hypersurfaces in $\L^{n+1}$, and thus we identify the tangent space of $M \subset\L^{n+1}$ at $p \in M$, denoted by $T_pM$, with a vector subspace of dimension $n$ in $\L^{n+1}$. In particular, by abuse of notation, if $\phi: U  \to M$, where $U$ is an open subset of  $\R^n$, is a local parametrization, we still use the symbol $\partial_i$ to denote the vector $\frac{\partial \phi} {\partial x_i}$. 

\begin{definition}\label{deftipology2}
Let $M \subset \L^{n+1}$ be a hypersurface. We say that $M$ is spacelike (resp. timelike, lightlike) if, for any $p \in M$, the vector subspace $T_pM$ is spacelike (resp. timelike, lightlike).  We say that $M$ is a non-degenerate hypersurface if $M$ is spacelike or timelike.
\end{definition}


\begin{definition}
A timelike vector $v \in \L^{n+1}$ is said to be future-oriented (resp. past-oriented) if $\l v, E_{n+1}\r <0$ (resp. $\l v, E_{n+1}\r >0)$, where $E_{n+1}:=(0,\ldots,0,1)$.
\end{definition}

We observe that for a spacelike (resp. timelike) surface $M$ and $p\in M$, we have the decomposition $\L^{n+1}=T_p M \oplus (T_p M)^\perp$, where $(T_p M)^\perp$ is a timelike (resp. spacelike) subspace of dimension 1 (see \cite{Lopez2014}).  A Gauss map is a differentiable map $N:M \to \L^{n+1}$ such that $|N(p)|=1$ and $N(p) \in (T_p M)^\perp$ for all $p \in M$. If $M$ is spacelike, the Gauss map pointing to the future is a map $N:M \to \H^{n}$. 


We recall now a result which is simple but crucial because it marks a relevant difference between the euclidean geometry and the geometry of Lorentz-Minkowski spaces.

\begin{proposition}\label{nonexistclosedspacelike}
Let $M \subset \L^{n+1}$ be a compact spacelike, timelike or lightlike hypersurface. Then $\partial M \neq \emptyset$.
\end{proposition}

\begin{proof}
Assume that $\partial M = \emptyset $ and that $M$ is spacelike (resp. timelike or lightlike). Let $a \in \L^{n+1}$ be a spacelike (resp. timelike) vector. Since $M$ is compact, there exists a minimum (or a maximum) $p_0 \in M$ for the function $f(p)=\l p, a\r$. Since $\partial M = \emptyset$, then $p_0$ is a critical point of the function $f$ and thus $\l v, a\r=0$ for all $v \in T_pM$. Hence $a \in (T_p M)^\perp$, but this gives a contradiction because $(T_p M)^\perp$ is timelike (resp. spacelike or lightlike).
\end{proof}

In other words, the previous result tells us that a closed hypersurface (i.e. compact without boundary) must be degenerate (see Definition \ref{deftipology2}). Therefore closed surfaces are not relevant in the Lorentz-Minkowski space, and this is deeply in contrast to the euclidean geometry.  For the sake of completeness, we also point out that Proposition \ref{nonexistclosedspacelike}, as well the previous definitions, can be extended to general hypersurfaces, see e.g. \cite[Sect. 3]{Lopez2014}.

\begin{definition}\label{defhyperbolicspace}
The hyperbolic space of center $p_0 \in \L^{n+1}$ and radius $r>0$ is the hypersurface defined by
$$ \H^n(p_0,r):=\{p \in \L^{n+1}; \ \l p- p_0, p-p_0\r = - r^2, \ \l p-p_0, E_{n+1} \r <0\},$$
where $E_{n+1}=(0,\ldots,0,1)$.
\end{definition}

From the euclidean point of view, this hypersurface is the ``upper sheet'' of a hyperboloid of two sheets.

\begin{remark}\label{hyperbolicspacelike}
The hyperbolic space is a spacelike hypersurface (see \cite{Lopez2014,Spivak}). In fact, let $v \in T_p\H^n(p_0,r)$ and let $\sigma=\sigma(s)$ be a curve in $\H^n(p_0,r)$  such that $\sigma^\prime(0)=v$. Then, differentiating with respect to $s$ the relation $\l \sigma(s) - p_0, \sigma(s) - p_0 \r = - r^2$ at $s=0$, we obtain $\l v, p-p_0 \r =0$. This implies that $T_p\H^n(p_0,r)= Span\{p-p_0\}^\perp$. Since $p-p_0$ is a timelike vector, it follows that $\H^n(p_0,r)$ is a spacelike hypersurface. Moreover $N(p)=(p-p_0)/r$ is a Gauss map. 
\end{remark}

When $p_0$ is the origin of $\L^{n+1}$, and $r=1$, the hyperbolic space is denoted by $\H^n$, that is 
$$\H^n:=\{(x_1,\ldots,x_{n+1}) \in \L^{n+1}; \ x_1^2+\ldots+x_n^2-x_{n+1}^2 = - 1, \ x_{n+1}>0\}.$$
In view of the previous remark, for any $p \in \H^n$, the induced metric on $T_p\H^n$ is positive definite, and hence the tensor $g=dx_1 \otimes dx_1+ \ldots +dx_n \otimes dx_n - dx_{n+1} \otimes dx_{n+1}$ is a Riemannian metric for $\H^n$. 
Another model for $\H^n$ is the Poincar\'e model in the unit disk $\B^n:=\{y \in \R^n; \|y\|_n<1\}$, where $\|\cdot\|_n$ is the euclidean norm in $\R^n$. The hyperbolic metric in $\B^n$ is defined by
$$\tilde g=\frac{4}{(1-\|y\|_n^2)^2} \sum_{i=1}^ndy_i\otimes dy_i,$$
which is conformally equivalent to the flat metric in $\B^n$. 
The isometry between $(\H^n, g)$ and $(\B^n, \tilde g)$ is given by the map $F:\H^n \to \B^n$ defined by
\begin{equation}\label{hyperbolicstereoproj}
F(x):= x_0 - \frac{2(x-x_0)}{\l x-x_0, x-x_0\r}=\left(\frac{x_1}{1+x_{n+1}},\ldots,\frac{x_n}{1+x_{n+1}}\right),
\end{equation}
where $x_0=(0,\ldots,0,-1) \in \R^{n+1}$ (see \cite[Proposition 3.5]{Lee}). The map $F$ is also known as hyperbolic stereographic projection, and from a geometrical point of view, $F$ sends a point $x \in \H^n$ to the intersection between the line joining $x$ and $x_0$ with the hyperplane $\{y \in \R^{n+1}; \ y_{n+1}=0\}$.



We conclude this section by recalling a variant of the Leray-Schauder fixed point theorem which will be used in the proof of Theorem \ref{mainteoprob}.

\begin{theorem}[A. J. B. Potter \cite{Potter}]\label{LSPTheom}
Let $X$  be a locally convex linear Hausdorff topological space and $U$ a closed convex subset of $X$ such that the zero element of $X$ is contained in the interior of $U$. Let $T:[0, 1]\times U \to X$ be a continuous map such that $T([0, 1] \times U )$ is relatively compact in $X$. Assume that
\begin{itemize}
\item[a)] $T(t,x)\neq x$ for all $x \in \partial U$ and $t \in [0,1]$;
\item[b)]	$T({0}\times \partial U)\subset U$. 
\end{itemize}
Then, there is an element $\bar x$ of $U$ such that $\bar x = T(1,\bar x)$.
\end{theorem}

\section{Derivation of the equation}
Let $\Omega$ be a proper smooth bounded domain of the hyperbolic space $\H^n$. Let us denote by $\mathcal{T}(\Omega)$ the space of tangent vector fields to $\Omega$ and denote by $\nabla^0$ the Levi-Civita connection of $\L^{n+1}$. We recall that $\nabla^0$ coincides with the flat connection of $\R^{n+1}$, and we denote by $\nabla$ the induced Levi-Civita connection on $\Omega$. 
Let $u$ be a smooth function defined on $\Omega$. We denote, respectively, by $du$, $\nabla u$ the differential and the gradient of $u$, which is the only vector field on $\Omega$ such that $$du (X) = \l X, \nabla u\r,  \ \ \hbox{for any} \ X \in \mathcal{T}(\Omega).$$

 The second covariant derivative of $u$ is defined as
$$ \nabla_{X,Y} u:= \nabla_X \nabla_Y u - \nabla_X Y (u)= \nabla_X \nabla_Y u - \nabla_{{\nabla_X Y}} u, \ \ \hbox{for any} \ X,Y \in \mathcal{T}(\Omega),$$
and the Hessian of $u$, denoted by $\nabla^2 u$, is the symmetric 2-tensor given by
$$\nabla^2 u\ (X,Y):=\nabla_{X,Y} u, \ \ \hbox{for any} \ X,Y \in \mathcal{T}(\Omega).$$
The Laplacian of $u$, denoted by $\Delta u$, is the trace of the Hessian.\\

Let $\{e_1,\ldots,e_n\}$ be a local orthonormal frame field for $\Omega$ 
and let $\{\omega^1,\ldots,\omega^n\}$ be the dual coframe field, i.e. $\omega^i(e_j)=\delta_{ij}$, for any $i,j=1,\ldots,n$.
The connection forms $\omega_{ij}$'s are defined by

\begin{equation}\label{defconnform}
\omega_{ij}(X):=\l \nabla_X e_j , e_i\r,\ X \in \T(\Omega),
\end{equation}
 and thus we have
\begin{equation}\label{dercoveiej}
 \nabla_{e_i}e_j = \sum_{k=1}^n \omega_{kj}(e_i) e_k.
\end{equation}
We also recall that the connection forms are skew symmetric, i.e. $\omega_{ij}+ \omega_{ji}=0$, for any $i,j \in \{1,\ldots,n\}$. In terms of the dual coframe field the exterior derivative of $u$ (i.e. the differential) can be written as
$$ du=\sum_{i=1}^nu_i \omega^i,$$
where $u_i$ denotes the covariant derivative $\nabla_{e_i}u$. We will also use the notation $\nabla_i$ to denote $\nabla_{e_i}$.

For the second covariant derivatives, taking $X=e_i$, $Y=e_j$ and using \eqref{defconnform} we have 

\begin{equation}\label{covarderij}
\nabla_{e_i,e_j} u= \nabla_{e_i}u_j -\sum_{k=1}^n\omega_{kj}(e_i) u_k.
\end{equation}

From now on we will use the notation $u_{ij}$ to denote $\nabla_{e_i,e_j} u$. In particular the Hessian of $u$ can be written as $u_{ij} \omega_j  \otimes \omega_i$ and the Laplacian of $u$ as $ \Delta u= \sum_{i=1}^n u_{ii}$.

\begin{definition}
Let $A \subset \L^{n+1}$, we define the cone spanned by $A$ as the set
$$\mathcal{C}_A:=\{\rho q \in \L^{n+1}; \ q\in A,\ \rho >0\}.$$
\end{definition}

\begin{remark}\label{remarkframe}
Observe that setting $e_{n+1}(x):=x/|x|$, for $x \in \mathcal{C}_\Omega$, and extending the $e_i$'s as constant along radii, i.e. $e_i(x)=e_i(x/|x|)$, $x \in \mathcal{C}_\Omega$, for $i=1,\ldots,n$, we get that $\{e_1,\ldots,e_{n+1}\}$ is a local orthonormal frame field for $\mathcal{C}_\Omega$, where $e_{n+1}$ is the future oriented unit radial direction, i.e. $\l e_{n+1}, e_{n+1}\r =-1$, $\l e_{n+1},E_{n+1}\r<0$. We also observe that by direct computation we have $\nabla^0_{i} e_{n+1}=e_i$, for any $i=1,\ldots,n$. 
We remark that by definition $e_{n+1}(q)=q$ for any $q \in \Omega$, and by abuse of notation when writing $\nabla^0_w q$, where $w \in \mathcal{T}(\R^{n+1})$, it will be always understood that we are computing $\nabla^0_w e_{n+1}$ at $x=q$, as well as $\nabla^0_q w$ will stand for $\nabla^0_{e_{n+1}} w$. 
\end{remark}

In order to derive the equation of spacelike $H$-bumps one can argue as in \cite[Sect. 1]{TW} with minor adjustments. Indeed, we only need to take into account the changes due to the bilinear form $\l\cdot,\cdot\r$, and the definition of mean curvature for spacelike hypersurfaces \cite[Sect. 3.2]{Lopez2014}). For the sake of completeness we derive the equation following the scheme of \cite[Sect. 2]{Lo03}. 

Let $u \in C^0(\overline{\Omega}) \cap C^2(\Omega)$, let $\Sigma$ be the associated radial graph and let $\mathcal{Y}:\Omega \to \R^{n+1}$ the map defined as $\mathcal{Y}(q):=e^{u(q)}q$.
From Remark \ref{remarkframe} it holds that $\nabla^0_i q=e_i$ and thus
\begin{equation}\label{eq1covder}
\nabla^0_i \mathcal{Y} = \nabla^0_i (e^u q)= e^u u_i q + e^u e_i.
\end{equation}
Therefore a local basis for $T_{\mathcal{Y}(q)}\Sigma$ is given by
$$ E_i(q)= e^u (e_i + u_i q), \ \ i=1,\ldots,n,$$
and the components of the metric are 
$$g_{ij}=\l E_i, E_j\r=e^{2u} (\l e_i, e_j \r + u_i u_j \l q, q\r )=e^{2u} (\delta_{ij} - u_iu_j).$$
Since we look for a spacelike hypersurface we must have $|\nabla u|^2 < 1$, and by elementary computations we see that the inverse matrix $(g^{ij})$ is given by
\begin{equation}\label{invmatrixmetric}
g^{ij}=e^{-2u}\left(\delta_{ij}+\frac{u_iu_j}{1-|\nabla u|^2}\right).
\end{equation}
For the Gauss map we have
$$N(\mathcal{Y}(q))= \frac{q+ \sum_{i=1}^n u_k e_k }{(1-|\nabla u|^2)^{1/2}}.$$
Indeed it is elementary to verify that $\l N(\mathcal{Y}(q)), E_i\r=0$ for any $i=1,\ldots,n$ and
$$ \l N(\mathcal{Y}(q)), N(\mathcal{Y}(q))\r = \frac{-1 + |\nabla u|^2}{1-|\nabla u|^2}=-1.$$
Moreover, as $u=0$ on $\partial\Omega$, there exists $q_1 \in \Omega$ such that $\nabla u(q_1)=0$ and by definition $N(\mathcal{Y}(q_1))=q_1$ and thus $\l N(\mathcal{Y}(q_1)), E_{n+1} \r < 0$. Therefore, since $N\circ \mathcal{Y} \in C^0(\overline\Omega,\R^{n+1})$ and $\Omega$ is connected, it follows that $N(\mathcal{Y}(\Omega)) \subset \H^n$, so that $N$ is future oriented.
The coefficients of the second fundamental form are given by
\begin{equation}\label{secondfundform}
\sigma_{ij}=\l N, \nabla^0_{i} \nabla^0_j \mathcal{Y}\r =\frac{e^{u}\left(-\delta_{ij} + u_iu_j -u_{ij}\right)}{(1-|\nabla u|^2)^{1/2}}.
\end{equation}

Indeed, recalling Remark \ref{remarkframe} and \eqref{eq1covder}, by direct computation we have
$$\nabla^0_{i} \left(\nabla^0_j \mathcal{Y} \right) = e^{u} \left( u_iu_j q + \nabla^0_i \nabla^0_j u\ q +  u_j e_i +  u_i e_j +  \nabla^0_i e_j\right).$$
Hence, by using the relations $\l e_i,e_j\r=\delta_{ij}$, $\l e_i, q\r=0$, and regrouping the terms, we deduce that
\begin{equation}\label{eqrev2}
 \l N, \nabla^0_{i} \nabla^0_j \mathcal{Y}\r= \frac{e^{u}}{(1-|\nabla u|^2)^{1/2}}\left(+u_iu_j - \nabla^0_i \nabla^0_j u + \l\nabla^0_i e_j, q\r  + \sum_{k=1}^n  u_k  \l\nabla^0_i e_j, e_k\r\right).
\end{equation}
Since $\l\nabla^0_i e_j, q\r = - \l e_j,  \nabla^0_i q\r = - \l e_j,  e_i\r=-\delta_{ij}$
and $$ \nabla^0_i \nabla^0_j u -\sum_{k=1}^n   u_k  \l\nabla^0_i e_j, e_k\r  =\nabla_i \nabla_j u -\sum_{k=1}^n  u_k  \l\nabla_i e_j, e_k\r=u_{ij},$$ then from \eqref{eqrev2} we finally get \eqref{secondfundform}.

At the end, from \cite[Definition 3.3]{Lopez2014}, the mean curvature of a spacelike hypersurface at $p=\mathcal{Y}(q) \in \Sigma$ is given by 
$$n H(\mathcal{Y}(q))= - \sum_{i,j=1}^n g^{ij} \sigma_{ij}.$$
Therefore, from \eqref{invmatrixmetric} and \eqref{secondfundform}, we deduce that $u$ must satisfy the following equation
$$ \sum_{i,j=1}^n \left((1-|\nabla u|^2)\delta_{ij} + u_i u_j)\right)u_{ij}= n (1-|\nabla u|^2) - n e^u (1-|\nabla u|^2)^{3/2}H(\mathcal{Y}(q)).$$

\section{A priori estimates}

Let $\e>0$ and let $\eta_\e \in C^\infty_0([0,+\infty))$ be such that $r\eta_\e \in C^\infty_0([0,+\infty))$, $r \mapsto \eta_\e(r)r$ is increasing in $(0,\frac{2}{\e})$ and decreasing in $(\frac{2}{\e}, +\infty)$. Assume moreover that 
$$\eta_\e(r)r=\left\{\begin{array}{cl}
  r & \hbox{for} \ r< 1-\e, \\[2mm]
  1-\frac{1}{2}\e  & \hbox{for} \  1-\frac{1}{2}\e<r<\frac{2}{\e},\\[2mm]
  0 & \hbox{for} \ r> \frac{3}{\e}.
\end{array}
\right.
$$
We define the regularized equation as
\begin{equation}\label{regularizedequation}
\begin{array}{lll}
&& \displaystyle \sum_{i,j=1}^n ((1- \eta^2_\e(|\nabla u|)|\nabla u|^2)\delta_{ij}+\eta_\e^2(|\nabla u|) u_iu_j) u_{ij}\\
&=&\displaystyle n (1-\eta_\e(|\nabla u|)^2 |\nabla u|^2) \left(1- \sqrt{1-\eta_\e(|\nabla u|)^2 |\nabla u|^2} e^u H(e^u q)\right).
\end{array}
\end{equation}
To simplify the notation we will write $ \eta^2_\e|\nabla u|^2$ instead of $\eta^2_\e(|\nabla u|)|\nabla u|^2$. The regularized Dirichlet problem for spacelike H-bumps is 
\begin{equation}\label{eqpseudoHbumps}
\begin{cases}
 \displaystyle  \sum_{i,j=1}^n  ((1- \eta^2_\e|\nabla u|^2)\delta_{ij}+\eta^2_\e u_iu_j) u_{ij}= n (1-\eta_\e^2 |\nabla u|^2) \left(1- \sqrt{1-\eta_\e^2 |\nabla u|^2} e^u H(e^u q)\right) & \hbox{in} \ \Omega,\\[6pt]
u=0 & \hbox{on}\  \partial \Omega.
\end{cases}
\end{equation}
We denote by $Q_{\e}$ the operator $$Q_{\e} (u):= \sum_{i,j=1}^n ((1- \eta^2_\e|\nabla u|^2)\delta_{ij}+\eta^2_\e u_iu_j) u_{ij} - n (1-\eta_\e^2 |\nabla u|^2) +n (1-\eta_\e^2 |\nabla u|^2)^{3/2} e^u H(e^u q).$$   
We claim that, in hyperbolic stereographic coordinates, the operator $Q_{\e}$ is uniformly elliptic. This is the content of the next lemma.
\begin{lemma}\label{lemmaunifelliptic}
For any $\e \in (0,1)$, the operator $Q_{\e}$, in hyperbolic stereographic coordinates, is uniformly elliptic with ellipticity constants depending only on $\e$ and $\Omega$.
\end{lemma}
\begin{proof}
Let $F:\H^n \to \B^n$ be the hyperbolic stereographic projection. By definition we have that $\Omega$ is mapped into a smooth proper domain $\Lambda=F(\Omega) \subset \subset \B^n$, $\phi=F^{-1}:\Lambda \to \Omega$ is a global parametrization, and there exist $c_1, c_2>0$ depending only on $\Omega$ such that 
\begin{equation}\label{limitazcost}
 c_1\leq \frac{(1-\|y\|_n^2)^2}{4} \leq c_2,\  \hbox{for all}\ y \in \Lambda.
\end{equation}
Let us set 
\begin{equation} \label{deflambda}
\lambda(y):=\frac{2}{1-\|y\|^2_n}, \ \ y \in \B^n.
\end{equation}
We recall that $F$ is an isometry and the hyperbolic metric in $\B^n$ is $\tilde g= \lambda^2 \sum_{i=1}^ndy_i\otimes dy_i$ (see Sect. 2). In particular, $\l \partial_i , \partial_j \r=\delta_{ij} \lambda^2$, where $\partial_i$ denotes the vector $\frac{\partial \phi} {\partial y_i}$, and the Christoffel symbols of the hyperbolic Levi-Civita connection are given by  
\begin{equation}\label{eqChristSymbrev}
\Gamma_{ij}^k=\frac{\lambda_i}{\lambda}\delta_{jk}+\frac{\lambda_j}{\lambda}\delta_{ik}-\sum_{l=1}^n\delta_{kl}\frac{\lambda_l}{\lambda}\delta_{ij},
\end{equation}
where $\lambda_i=\frac{\partial \lambda}{\partial y_i}$.
In local coordinates the gradient is given by
\begin{equation}\label{formulagradrev}
\nabla u=\sum_{i,j=1}^n \tilde g^{ij} \frac{\partial \tilde u}{\partial y_i} \partial_j=  \lambda^{-2} \sum_{i=1}^n \frac{\partial \tilde u}{\partial y_i} \partial_i,
\end{equation}
and thus
\begin{equation}\label{formulatrasfgrad} 
\begin{array}{lll}
\displaystyle |\nabla u|^2&=&\displaystyle  \lambda^{-2}\|\nabla_0 \tilde u\|_n^2,
\end{array}
\end{equation}
where $\tilde u= u \circ F^{-1}$ and $\nabla_0 \tilde u$ is the gradient of $\tilde u$ with respect to the euclidean flat metric. Using the well known expression for the Hessian and the Laplacian in local coordinates we have

\begin{equation}\label{eqhessianlocal}
\nabla^2 u (\partial_i,\partial_j) = \frac{\partial^2\ \tilde u}{\partial y_j \partial y_i} - \sum_{k=1}^n \Gamma_{ji}^k \frac{\partial \tilde u}{\partial y_k},
\end{equation}
 
\begin{equation}\label{eqlaplacianhypcord}
\Delta u = \sum_{i,j=1}^n {\tilde g}^{ij} \left(\frac{\partial^2 \tilde u}{\partial y_i \partial y_j} - \sum_{k=1}^n \Gamma_{ij}^k \frac{\partial \tilde u}{\partial y_k}\right) =  \lambda^{-2} \sum_{i=1}^n \left(\frac{\partial^2 \tilde u}{\partial y_i^2} - \sum_{k=1}^n \Gamma_{ii}^k \frac{\partial \tilde u}{\partial y_k}\right).
\end{equation}
By using the previous identities and \eqref{eqChristSymbrev} we infer that
\begin{equation}\label{partediordine2}
\sum_{i,j=1}^n u_iu_j u_{ij} = \lambda^{-4} \sum_{h,k=1}^n  \frac{\partial \tilde u}{\partial y_h} \frac{\partial \tilde u}{\partial y_k} \frac{\partial^2 \tilde u}{\partial y_h \partial y_k} + \Phi,
\end{equation}
where $\Phi$ is a term which does not involve second order partial derivatives. From \eqref{deflambda}, \eqref{formulatrasfgrad}, \eqref{eqlaplacianhypcord} and \eqref{partediordine2} we deduce that the principal part of the operator $Q_\e$, in hyperbolic stereographic coordinates, is
$$\lambda^{-2} \left[\sum_{i,j=1}^n \left(1-\eta_\e^2  \lambda^{-2} \|\nabla_0  \tilde u\|_n^2\right) \delta_{ij} \frac{\partial^2 \tilde u}{\partial y_i \partial y_j} + \eta_\e^2 \lambda^{-2} \frac{\partial \tilde u}{\partial y_i} \frac{\partial \tilde u}{\partial y_j} \frac{\partial^2 \tilde u}{\partial y_i \partial y_j}\right],$$
where $\eta_\e=\eta_\e(\lambda^{-1}\|\nabla_0  \tilde u\|_n)$.
For any $i,j=1,\ldots,n$, we define, for $y \in \Lambda$, $p=(p_1,\ldots,p_n) \in \R^n$,
\begin{equation}\label{coeffelliptcity}
\tilde a^{ij}_\e(y,p):=\lambda^{-2}\left[ \left(1-\eta_\e^2(\lambda^{-1} \|p\|_n) \lambda^{-2} \|p\|_n^2\right) \delta_{ij} + \eta_\e^2(\lambda^{-1} \|p\|_n)\lambda^{-2}p_ip_j\right].
\end{equation}
Now, for any $\xi=(\xi_1,\ldots,\xi_n) \in \R^n$, $y \in \Lambda$, $p \in \R^n$ we claim that
\begin{equation}\label{eq1unifelliptc}
c_2\|\xi\|_n^2\geq \sum_{i,j=1}^n \tilde a^{ij}_\e(y,p) \xi_i \xi_j \geq \frac{1}{2}\e c_1\|\xi\|_n^2,
\end{equation}
where the constants $c_1,c_2$ are given by \eqref{limitazcost}.
Indeed by the definition of $\eta_\e$  for any $y \in \Lambda$, $p \in \R^n$  it holds 
$$0\leq \eta_\e^2(\lambda^{-1} \|p\|_n) \lambda^{-2} \|p\|_n^2 \leq (1-\e/2)^2$$
and thus
\begin{eqnarray*}
 \sum_{i,j=1}^n \tilde a^{ij}_\e(y,p) \xi_i \xi_j &=&\lambda^{-2} \left[ \left(1-\eta_\e^2 \lambda^{-2} \|p\|_n^2\right) \|\xi\|_n^2 + \eta_\e^2\lambda^{-2}\left(\sum_{i=1}^n\ p_i \xi_i\right)  \left(\sum_{j=1}^np_j\xi_j\right)\right]\\
&=& \lambda^{-2} \left[ \left(1-\eta_\e^2 \lambda^{-2}\|p\|_n^2\right) \|\xi\|_n^2 + \eta_\e^2\lambda^{-2} (p,\xi)_{n}^2\right]\\
&\geq&  \lambda^{-2} \left(1-\eta_\e^2 \lambda^{-2}\|p\|_n^2\right) \|\xi\|_n^2\geq  \lambda^{-2} (1-(1-\e/2)^2)\|\xi\|_n^2 \geq \frac{1}{2}c_1\e \|\xi\|_n^2,
\end{eqnarray*}
where $(\cdot,\cdot)_{n}$ denotes the euclidean scalar product in $\R^n$. The proof of the other inequality in \eqref{eq1unifelliptc} is similar and we omit the details. The proof is then complete.
\end{proof}
For $t \in [0,1]$, we define the operator
$$Q^t_{\e} (u):=\sum_{i,j=1}^n((1- \eta^2_\e|\nabla u|^2)\delta_{ij}+\eta^2_\e u_iu_j) u_{ij} - n t (1-\eta_\e^2 |\nabla u|^2) +n t (1-\eta_\e^2 |\nabla u|^2)^{3/2} e^u H(e^u q),$$   
For $u$ such that $|\nabla u|_{\infty,\Omega}<1$, we also define the operator
$\mathcal{Q}^t(u)$ as
\begin{equation}\label{operatorcomp}
\mathcal{Q}^t(u):=\sum_{i,j=1}^n((1- |\nabla u|^2)\delta_{ij}+u_iu_j) u_{ij} - nt (1- |\nabla u|^2) +nt (1- |\nabla u|^2)^{3/2} e^u H(e^u q).
\end{equation}  

\begin{remark}\label{remoperator}
By definition, for any fixed $\e \in (0,1)$, if $u$ is such that $|\nabla u|_{\infty,\Omega} \leq 1-\e$ we have  $Q^t_{\e} (u)=\mathcal{Q}^t(u)$, for any $t \in [0,1]$. Moreover, in view of Lemma  \ref{lemmaunifelliptic} and since the principal parts of $Q_\e^t$, $\mathcal{Q}^t$ are independent on $t$, they are uniformly elliptic even with respect to $t$, when passing to hyperbolic stereographic coordinates.\\ 
\end{remark}
\begin{remark}\label{remarkstrutturaoperatoretrasf}
As seen in the proof of Lemma \ref{lemmaunifelliptic} we can write an explicit expression of the operator $Q^t_\e$ in hyperbolic stereographic coordinates defined in the whole $\Lambda=F(\Omega)$. For our purposes we just observe that the transformed operator is of the form
$$\tilde{Q}^t_\e(\tilde u)=\sum_{i,j=1}^n \tilde a_\e^{ij}(y,\nabla_0 \tilde u) \tilde u_{ij} + \tilde b_{\e,t}(y,u,\nabla_0 \tilde u),$$  
where $\tilde a_\e^{ij}=\tilde a_\e(y,p):\Lambda \times \R^n \to \R$ is given by \eqref{coeffelliptcity}, and $\tilde b_{\e,t}:\Lambda \times \R \times \R^n \to \R$ is given by
\begin{eqnarray*}
\tilde b_{\e,t}(y,z,p)&:=&-\left(1-\eta_\e^2  \lambda^{-2} \|p\|_n^2\right)\sum_{k=1}^n G_k(y) p_k - \eta_\e^2 \lambda^{-2}\sum_{h,k,r=1}^n G_{hkr}(y) p_h p_k p_r\\
&&- nt \left(1-\eta_\e^2 \lambda^{-2}\|p\|_n^2\right) + nt \left(1-\eta_\e^2 \lambda^{-2} \|p\|_n^2\right)^{3/2}e^zH(e^zF(y)),
\end{eqnarray*}
where $\lambda$ is defined in \eqref{deflambda}, $G_k$,$G_{hkr}$ are smooth functions defined in $\Lambda$, $h,k,r \in \{1,\ldots,n\}$. We point out that $\tilde a_\e^{ij}$ does not depend on $z$ and $\tilde a_\e^{ij}=O(1)$, $\tilde b_{\e,t}=O(\|p\|_n)$, as $\|p\|_n\to + \infty$, uniformly for $y \in \Lambda$, and $z$ in compact subsets of $\R$. In particular, according to the notations of \cite{GT}, setting  $\mathcal{E}:=\sum_{ij=1}^n \tilde a_\e^{ij} p_ip_j$, we have that $\mathcal{E}$ does not depend on $z$ and $\mathcal{E}=O(\|p\|_n^2)$ as $\|p\|_n \to + \infty$, uniformly for $y \in \Lambda$. 

These properties will be useful in the sequel. In addition, since $(\tilde a_\e^{ij})$ is symmetric and positive definite, when applying the results of \cite[Sect. 15]{GT}, it will be understood that we take $(\tilde a_\e^{ij})^*=\tilde a_\e^{ij}$ and $c_i=0$ (see \cite[(15.23)]{GT}). 
\end{remark}

We define now the class of admissible domains.
\begin{definition}\label{classadmissible}
Let $\Omega$ be a bounded domain of $\H^n$ and let $H \in C^0(\mathcal{C}_{\overline{\Omega}})$. We say that $(\Omega,H)$  is admissible if there exists a constant $\theta \in (0,1)$ such that for any $q_0 \in \partial \Omega$ and for any $t \in [0,1]$, there exist two functions $\varphi_1, \varphi_2 \in C^{2}(\overline\Omega)$ satisfying 
\begin{itemize}
\item[(i)] $\sup_\Omega |\nabla \varphi_i| \leq 1-\theta$,  for $i=1,2$,
\item[(ii)] $\varphi_1(q_0) =0$ and $\varphi_1(q_0) \leq 0$ on $\partial \Omega$,
\item[(iii)]  $\varphi_2(q_0) =0$ and $\varphi_2(q_0) \geq 0$ on $\partial \Omega$,
\item[(iv)]  $\mathcal{Q}^t (\varphi_1) \geq 0  $, $\mathcal{Q}^t (\varphi_2) \leq 0 $ in $\Omega$.
\end{itemize}
We denote by $\mathcal{A}$ the set of admissible couples $(\Omega,H)$. Given $\theta \in (0,1)$, and given $\Omega$ and $H$ as above, we say that $(\Omega,H,\theta)$ is admissible if $(\Omega,H)$ is admissible with constant $\theta$.
\end{definition}

 \begin{remark}\label{remadmissible}
We observe that $\mathcal{A}\neq \emptyset$. In fact for any given domain $\Omega\subset \H^n$, for any fixed $m>0$, the function $H(x)=\frac{1}{|x|^m}$, $x \in \mathcal{C}_{\overline\Omega}$ is such that $(\Omega,H) \in \mathcal{A}$. In fact it is easy to see that $Q_t(0)=0$ for any $t \in [0,1]$, so that the functions $\varphi_1=0$, $\varphi_2=0$ verify (i)-(iv) for any $\theta \in (0,1)$. More in general, for any domain $\Omega$, for any function $H \in C^0(\mathcal{C}_{\overline{\Omega}})$ such that $H{\big |_\Omega}=1$, we have that $(\Omega,H) \in \mathcal{A}$, and $(\Omega,H,\theta)$ is admissible for any $\theta \in (0,1)$. 
\end{remark}
This condition of admissibility is very general. If we impose some regularity on $\partial\Omega$, and if we assume that $H$ is positive, smooth and not increasing along radii, then every couple $(\Omega, H)$ is admissible. This is the content of the next result. We introduce first the following definition.

\begin{definition}\label{defextgeodesic}
Let $\Omega$ be a bounded domain of $\H^n$. We say that $\Omega$ satisfies a uniform exterior geodesic ball condition if there exist $\sigma>0$ 
and a map  $\Xi:\partial\Omega \to \H^n$ of class $C^2$ such that for any $q_0 \in \partial \Omega$ there exists a geodesic ball in $\H^n$ of radius $\sigma$ centered at $\xi=\Xi(q_0) \in \H^n \setminus \overline \Omega$, and denoted by $B_\sigma(\xi)$, such that $q_0 \in \partial B_\sigma(\xi)$ and $B_\sigma(\xi) \subset \H^n \setminus \overline\Omega$.
\end{definition}

\begin{proposition}\label{propadmissextbgball}
Let $\Omega$ be a bounded domain of $\H^n$ satisfying a uniform exterior geodesic ball condition. Let $H \in C^1(\mathcal{C}_{\overline\Omega})$ be such that $H>0$ and $\frac{\partial}{\partial \lambda} \left(\lambda H(\lambda q) \right)\leq 0,\ \hbox{for all}\ q \in \overline\Omega, \lambda>0$. Then $(\Omega,H)$ is admissible
\end{proposition}
\begin{proof}
Let $dist_{\H^n}(\cdot,\cdot)$ be the geodesic distance in $\H^n$. Let $\sigma>0$ be the number given by Definition \ref{defextgeodesic} for $\Omega$. In particular, by definition, it follows that for any $q_0 \in \partial \Omega$ there exists $\xi=\xi(q_0) \notin \overline\Omega$ such that $dist_{\H^n}(\xi,\partial \Omega)=dist_{\H^n}(\xi,q_0)=\sigma$. 

Let $q_0 \in \partial \Omega$, $t \in [0,1]$ and let $\xi=\xi(q_0)$ satisfying the above properties. Since every geodesic ball of $\H^n$ is geodesically convex (see \cite[Sect 2.5]{Papadopoulos}) we can take $R>0$ sufficiently large so that $\overline\Omega$ is contained in the geodesically convex ball $B_R(\xi)$. We observe that since $\Omega$ is bounded and $dist_{\H^n}(\xi, \partial\Omega)=\sigma$, up to a new choice of a larger $R$, we can assume that $R$ is uniform with respect to the choice of $q_0 \in \partial\Omega$. 

Arguing as in proof of \cite[Theorem 2.1]{Gerhardt}, we set $|||q|||:=dist_{\H^n}(q,\xi)$ to denote the geodesic distance from $\xi$ and we define
$$\delta^+(q):= \int_{|||q_0|||}^{|||q|||} (1+\gamma(s))^{-1/2} \ ds,$$ 
$$\gamma(s):=\alpha e^{\beta s},$$
where
$\alpha$, $\beta$ are positive constant to be determined later. By construction it holds that $\delta^+ \in C^2(\overline{B_R(\xi)}\setminus \{\xi\})$, $\delta^+ \in C^2(\overline\Omega)$, $\delta^+(q_0)=0$ and $\delta^+\geq0$ in $\overline\Omega$ because of the exterior ball condition.

Let us consider the following operator (which is the divergence form of $-\mathcal{Q}^t$)
$$\mathcal{Q}_{\mathrm{div}}^t(u):=-\mathrm{div}_{\H^n}\left(\frac{\nabla u}{\sqrt{1-|\nabla u|^{2}}}\right)+\frac{nt}{\sqrt{1-|\nabla u|^{2}}}-nte^{u}H(e^{u}q).$$
We set 
$$A(u):=-\mathrm{div}_{\H^n}\left(\frac{\nabla u}{\sqrt{1-|\nabla u|^{2}}}\right),\ \ \nu(u):=\frac{1}{\sqrt{1-|\nabla u|^{2}}}.$$

We observe that $|\nabla |||q||||=1$ for any $q \in {B_R(\xi)}\setminus\{\xi\}$. This property is known for general manifolds when $R$ is sufficiently small so that $B_R(\xi)$ is contained in a normal neighborhood of $\xi$ (see \cite[Corollary 6.9 and 6.11]{Lee}). In our case, as a consequence of the Cartan-Hadamard Theorem, since $\H^n$ has negative sectional curvature then it admits global normal coordinates and we are done. 

Therefore, since the covariant derivatives of $\delta^+$ are given by $(\delta^+)_i=(1+\gamma)^{-1/2}|||q|||_i$, we obtain that for any $q \in \overline\Omega$
$$|\nabla \delta^+|=(1+\gamma)^{-1/2} <1,$$
and 
$$\nu(\delta^+)=\gamma^{-1/2}(1+\gamma)^{1/2}.$$
In addition, by direct computation (see \cite[(2.14)-(2.16)]{Gerhardt}) it holds that
$$A(\delta^+)=(1+\gamma)^{-1/2}\left(\frac{\beta}{2}-\Delta|||q|||\right) \nu(\delta^+),$$
We observe that $\Delta|||q|||$ is smooth and bounded in compact subsets of $\overline{B_R}(\xi)\setminus \{\xi\}$ and it is singular as $q \to \xi$. Indeed (see \cite[(2.17)-(2.18)]{Gerhardt}) we have
\begin{equation}\label{eqLaplgeoddist}
-\Delta|||q|||=-\frac{n-1}{|||q|||}+\Psi,
\end{equation}
where $\Psi$ is a bounded term which is given, in normal coordinates centered at $\xi$, by
$$\Psi = - \sum_{i,j,k=1}^n g^{ij}\Gamma_{ij}^k |||q|||_k.$$  
In particular, in view of the uniform exterior geodesic ball condition, since $dist_{\H^n}(q,\xi(q_0))\geq \sigma$ for any $q \in \overline\Omega$, for any $q_0 \in \partial\Omega$, then from \eqref{eqLaplgeoddist} we infer that $\Delta|||q|||$ is bounded in $\overline\Omega$ by a constant depending only on $n, \sigma,\Omega,q_0$. In addition, by definition the map $q_0 \mapsto \xi$ is of class $C^2(\partial \Omega, \H^n)$ and thus, by compactness of $\partial \Omega$, it follows that $\Delta|||q|||$ is bounded by a constant depending only $n,\sigma,\Omega$.\\

Now, by the previous relations we have
\begin{eqnarray*}
A(\delta^+)+ tn\nu(\delta^+)&=& \left[(1+\gamma)^{-1/2}(\beta/2-\Delta|||x|||) + tn\right](\gamma^{-1/2}(1+\gamma)^{1/2})\\[6pt]
&\geq&\left(\frac{\beta}{2}-\Delta|||x|||\right)\gamma^{-1/2}= \left(\frac{\beta}{2}-\Delta|||x|||\right) \alpha^{-1/2} e^{-\frac{\beta s}{2}}.
\end{eqnarray*}
Setting $\bar H:= \max_{q \in \overline\Omega} H(q)>0$, we can choose $\beta$ sufficiently large so that $\frac{\beta}{2}-\Delta|||q|||>0$, for any $q \in \overline\Omega$. With this choice of $\beta$ we choose $\alpha$ sufficiently small so that  $\left(\frac{\beta}{2}-\Delta|||q|||\right) \alpha^{-1/2} e^{\frac{\beta s}{2}} \geq n \bar H$, for any $x \in \overline\Omega$, $s \in [|||q_0|||, \sup_{q \in \overline\Omega}|||q|||]$. Therefore, since $\delta^+\geq0$ in $\overline\Omega$, and in view of the monotonicity assumption on $H$, it follows that
\begin{eqnarray*}
A(\delta^+)+nt\nu(\delta^+)-nte^{\delta^+}H(e^{\delta^+}q) &\geq& A(\delta^+)+nt\nu(\delta^+)-nte^{0}H(e^{0}q)\\
& \geq& A(\delta^+)+nt\nu(\delta^+)-nt \bar H\\
& \geq& 0.
\end{eqnarray*}
Hence, $\mathcal{Q}_{\mathrm{div}}^t(\delta^+) \geq 0$ in $\overline\Omega$, which is equivalent to $\mathcal{Q}^t(\delta^+)\leq 0$ in $\overline\Omega$, and in addition by construction we have $\delta^+\geq 0$ on $\partial\Omega$, $\delta^+(x_0)=0$, $|\nabla \delta^+|=(1+\gamma)^{-1/2}\leq 1 - \theta_+$, for some number $\theta_+=\theta_+(\alpha, \beta) \in (0,1)$. 

As pointed out before, in view of the uniform exterior ball condition, $-\Delta|||q|||$ is uniformly bounded by a constant depending only on $n,\sigma,\Omega$, and by construction $\sup_{q \in \overline\Omega}|||q|||\leq R$. Therefore, the numbers $\alpha$, $\beta$ can be chosen in a uniform way with respect to the base point $q_0 \in \partial \Omega$ (and also with respect on $t \in [0,1]$). Hence, there exist $\theta_+ \in (0,1)$ such that for any $q_0 \in \partial \Omega$, $t \in[0,1]$ the function $\varphi_2:=\delta^+$ (which depends on the choice of $q_0$ but non on $t$) satisfies i)-iv) of Definition \ref{classadmissible} with $\theta=\theta^+$. For the other barrier is suffices to take $\varphi_1:=\delta^-$, where
 $$\delta^-:=- \int_{|||q_0|||}^{|||q|||} (1+\gamma(s))^{-1/2} \ ds,$$ 
and to argue as in the previous case. We limit to observe that in this case the choice of $\alpha$, $\beta$ has to be made in a different way but it is still uniform with respect to $q_0$, and $t$. 

In fact
\begin{eqnarray*}
A(\delta^-)+ tn \nu(\delta^-)&=& - \gamma^{-1/2}\left(\frac{\beta}{2}-\Delta|||x|||- tn(1+\gamma)^{1/2}\right)\\[6pt]
&=& - \alpha^{-1/2} e^{-\frac{\beta s}{2}}\left(\frac{\beta}{2}-\Delta|||x|||- tn(1+\alpha e^{\beta s})^{1/2}\right).
\end{eqnarray*}
Taking $\alpha=e^{- \beta \sup_{x \in \overline\Omega}|||x|||}$, it follows that $n(1+\alpha e^{\beta s})^{1/2}\leq \sqrt{2}n$ for any $\beta>0$, $s \in [|||x_0|||,\sup_{x \in \overline\Omega}|||x|||]$. With this choice of $\alpha$, we choose $\beta$ such that
$$\frac{\beta}{2}-\Delta|||x|||-2n \geq 0,$$
for $x \in \overline\Omega$. At the end, we have that $A(\delta^-)+ tn \nu(\delta^-) \leq 0,$ and thus being $H>0$ it holds that
$$A(\delta^-)+ tn \nu(\delta^-) - tn e^{\delta^-}H(e^{\delta^-}) \leq A(\delta^-)+ tn \nu(\delta^-) \leq 0 \ \hbox{in} \ \Omega.$$
As before we find a uniform $\theta_- \in (0,1)$ such that for any $q_0\in \partial \Omega$, $t \in [0,1]$, the function $\varphi_1:=\delta^-$ satisfies $|\nabla \varphi_1|\leq 1-\theta_-$ and ii)-iv) of Definition \ref{classadmissible}. At the end, choosing $\theta:=\min\{\theta_-,\theta_+\}$ we have that for any $q_0 \in \partial \Omega$, $t \in [0,1]$, the functions $\varphi_1$, $\varphi_2$ satisfy i)-iv) of Definition \ref{classadmissible}, and hence $(\Omega,H)$ is admissible. The proof is complete.
\end{proof}

\begin{remark}\label{rmk:thetadipendenceH}
It is important to note that in the previous proof the choice of $\theta$ depends only on $n,\sigma,\Omega$ and depends on $H$ just by the number $\bar H:= \max_{q \in \overline\Omega} H(q)>0$ because of the monotonicity assumption. In particular $\theta$ does not depend on the derivatives of $H$.
\end{remark}

 If $H \in C^{1}({C_{\overline{\Omega}}(r_1,r_2)})$ we define a canonical extension of $H$ to a mapping on the cone $\mathcal{C}_{\overline{\Omega}}$ in the following way: set
$$h_1(q):=\left[\frac{\partial}{\partial \rho} \rho H(\rho q)\right]_{\rho=r_1}, \ \ h_2(q):=\left[\frac{\partial}{\partial \rho} \rho H(\rho q)\right]_{\rho=r_2}$$
and
\begin{equation} \label{eq:extension}
\hat H(\rho q):=
\begin{cases}
\frac{r_1}{\rho} H(r_1q) + \left(1-\frac{r_1}{\rho}\right)h_1(q) & \hbox{for} \ \rho \in (0,r_1),\\
H(\rho q) &  \hbox{for} \ \rho \in [r_1,r_2],\\
\frac{r_2}{\rho} H(r_2q) + \left(1-\frac{r_2}{\rho}\right)h_2(q) & \hbox{for} \ \rho \in (r_2, +\infty).
\end{cases}
\end{equation}

\begin{remark}\label{remextensadmissible}
It is elementary to check that $\hat H \in C^1(\mathcal{C}_{\overline{\Omega}})$, and if $H$ satisfies $\frac{\partial}{\partial \lambda} \left(\lambda H(\lambda q) \right)\leq 0,\ \hbox{for all}\ q \in \overline\Omega, \lambda \in [r_1,r_2]$ it follows that
$$\frac{\partial}{\partial \lambda} \left(\lambda \hat H(\lambda q) \right)\leq 0,\ \hbox{for all}\ q \in \overline\Omega, \lambda >0.$$ 
Therefore, since $\hat H(x)=H(x)$ for $x \in \overline\Omega$, by Remark \ref{rmk:thetadipendenceH} if $\Omega$ satisfies the hypotheses of Proposition \ref{propadmissextbgball} and $H$ is positive, it follows that $(\hat H,\Omega)$ is admissible with constant which does not depend on the choice of $r_1,r_2$, and the derivatives of $\hat H$. 
\end{remark}

In view of the previous remark, the following definition makes sense:
\begin{definition}\label{classadmissibleext}
Let $\Omega$ be a bounded domain of $\H^n$, let $0<r_1 \leq 1 \leq r_2$ and $H \in C^{1}({C_{\overline{\Omega}}(r_1,r_2)})$. We say that $(\Omega, H)$ is admissible if $(\Omega, \hat H)$ is admissible, where $\hat H$ is the extension of $H$ defined in \eqref{eq:extension}, and for $\theta \in (0,1)$ we say 
that $(\Omega, H,\theta)$ is admissible if $(\Omega, \hat H,\theta)$ is admissible. 
\end{definition}

Now we have all the tools to prove the a priori estimates. Let us fix some notation: let $k \in \mathbb{N}$, $\alpha \in (0,1)$ we consider the subspaces $C^{k,\alpha}_{0}(\overline\Omega):=\{u\in C^{k,\alpha}(\overline\Omega)~|~u|_{\partial\Omega}=0\}$, $C^{k}_{0}(\overline\Omega):=\{u\in C^{k}(\overline\Omega)~|~u|_{\partial\Omega}=0\}$, endowed, respectively, with the usual norms $|\cdot|_{k,\alpha}$, $|\cdot|_{k}$. We point out that $C^{k,\alpha}_{0}(\overline\Omega)$, $C^{k}_{0}(\overline\Omega)$ are closed subspaces of Banach spaces and thus they are Banach too. When needed we will specify also the domain in the norms, otherwise it will be understood that the domain is $\Omega$, moreover for the $C^{0}(\overline\Omega)$-norm we will use the notations $|\cdot|_\infty$, $|\cdot|_{\infty,\Omega}$, and $\|\cdot\|_{\infty}$, $\|\cdot\|_{\infty,\Omega}$ when working in the euclidean setting.

 We define $\hat Q^t_\e$ as the operator obtained from $Q^t_\e$ by replacing $H$ with its extension $\hat H$, and in the class of functions satisfying $|\nabla u|_{\infty,\Omega}<1$ we define $\hat{\mathcal{Q}}^t$ as
\begin{equation}\label{operatorcomp2}
\hat{\mathcal{Q}}^t(u):=\sum_{i,j=1}^n((1- |\nabla u|^2)\delta_{ij}+u_iu_j) u_{ij} - nt (1- |\nabla u|^2) +nt (1- |\nabla u|^2)^{3/2} e^u \hat H(e^u q).
\end{equation} 
  
 In order to simplify the notation we set $L_{\e,u} u:= \sum_{i,j=1}^n ((1- \eta^2_\e|\nabla u|^2)\delta_{ij}+\eta^2_\e u_iu_j) u_{ij}$. The first result we prove is about a priori $C^0$ estimate for solutions of $\hat Q^t_\e(u)=0$.
 
\begin{lemma}[\textbf{A priori $C^0$ estimates}] \label{lemmaczeroestimates}
Let $\Omega$ be a bounded domain and let $r_1 \neq r_2$ such that $0<r_1\leq 1\leq r_2$. Assume that $H \in C^1(\mathcal{C}_{\overline{\Omega}}(r_1,r_2))$ satisfies
\begin{equation}\label{assumptHstrong}
H(r_1 q) > r_1^{-1} \ \hbox{and} \ H(r_2 q) < r_2^{-1} \ \hbox{for any} \ q \in \overline{\Omega},
\end{equation}
and 
\begin{equation}\label{assumptmonotonicity}
 \frac{\partial}{\partial \lambda} \left(\lambda H(\lambda q) \right)\leq 0,\ \hbox{for all}\ q \in \overline\Omega, \lambda \in [r_1,r_2].
\end{equation} 
For $\e \in(0,1)$, for every $t \in [0,1]$, if $u \in C^2_0(\overline\Omega)$ is a solution of $\hat Q^t_\e(u)=0$ then $$\log r_1 \leq u(q) \leq \log r_2, \ \hbox{for every}\ q \in \overline \Omega.$$ 
\end{lemma}

\begin{proof}
Let us observe that since we are assuming \eqref{assumptHstrong}, \eqref{assumptmonotonicity}, it holds that 
\begin{equation}
\label{eq:strict-barriers2}
\hat H(x)>|x|^{-1}\textrm{ if $|x|\le r_{1},~x\in\mathcal{C}_{\overline\Omega}$}\quad\textrm{and}\quad \hat H(x)<|x|^{-1}\textrm{ if $|x|\ge r_{2},~x\in\mathcal{C}_{\overline\Omega}$}.
\end{equation}

Let $\e \in (0,1)$, let $t \in [0,1]$ and let $u \in C^2_0(\overline\Omega)$ such that $\hat Q^t_\e(u)=0$. By definition $u$ is a classical solution of the Dirichlet problem

\begin{equation}\label{sistemaaux2}
\begin{cases}
L_{\e,u} u = n t (1-\eta_\e^2 |\nabla u|^2) \left(1- \sqrt{1-\eta_\e^2 |\nabla u|^2} e^u  \hat H(e^u q)\right) & \hbox{in} \ \Omega,\\
u=0 & \hbox{in} \ \partial\Omega.
\end{cases}
\end{equation}

Let $q_0 \in \overline \Omega$ such that $u(q_0)= \max_{\overline \Omega} u$. Assume by contradiction that $u(q_0) > \log r_2$. Then $q_0 \in \Omega$ because $r_2\geq 1$ and $u=0$ on $\partial \Omega$. Hence $\nabla u (q_0)=0$, and being $q_0$ a maximum point, it holds that $\Delta u (q_0) \leq 0$, and by definition of $L_{u, \e}$  this reads as
$$L_{u, \e} u \leq 0.$$
On the other hand it must be $t>0$ because otherwise if $t=0$ then $u\equiv 0$. Moreover
$$L_{u, \e} u (q_0)=nt e^{u(q_0)} \left(\frac{1}{e^{u(q_0)}}- \hat H(e^{u(q_0)}q_0)\right)>0,$$
because $\hat H(x) < |x|^{-1}$ as $|x|>r_2$. Thus we reach a contradiction. The same argument holds to show that $\min_{\overline \Omega} u \geq \log r_1.$
\end{proof}

\begin{lemma}[\textbf{A priori $C^{1,\alpha}$ estimates}]\label{lemapriorigradestimates}
Let $\e \in (0,1)$ and let $\Omega$ be a bounded domain of class $C^{2}$. Assume that $H$ satisfies \eqref{assumptHstrong}, \eqref{assumptmonotonicity}. Then, there exist  two positive constants $M$, $C$ and $\alpha_0 \in (0,1)$ such that for all $t \in [0,1]$
if $u \in C^2_0(\overline{\Omega})$ is such that $|\nabla u|_{\infty,\partial\Omega} \leq 1-\e$ and is a solution of the equation $\hat Q^t_\e(u)=0$, then
$$|\nabla u|_{\infty,\Omega} \leq M $$
$$|u|_{1,\alpha_0} \leq C.$$
\end{lemma}

\begin{proof}
Let us fix $\e \in (0,1)$, let $t \in [0,1]$ and let $u=u_t$ be a solution of $\hat Q^t_\e(u)=0$. From Lemma \ref{lemmaczeroestimates} we have that $\log r_1 \leq u \leq \log r_2$ and thus by definition $u$ also verifies $Q_\e^t(u)=0$. Therefore, from now on we can work just with the operator $Q_\e^t$.

 Let us set $b_{\e,t}(q,u,\nabla u):=n t (1-\eta_\e^2 |\nabla u|^2) \left(1- \sqrt{1-\eta_\e^2 |\nabla u|^2} e^u H(e^u q)\right)$.
In view of Remark \ref{remoperator}, passing to hyperbolic stereographic coordinates, the operator $Q_\e^t$ is uniformly elliptic with constants independent on $t$, moreover, by definition and thanks to Lemma \ref{lemmaczeroestimates} the term $b_{\e,t}(q,u,\nabla u)$ is uniformly bounded with respect to $t$. 

Now there are only two possibilities: there exists a constant $M$ independent on $t$ such that $|\nabla u|_{\infty,\Omega}\leq M$ for all $t \in [0,1]$ or there exists a subsequence $(t_k) \subset [0,1]$, such that $|\nabla u_{t_k}|_{\infty,\Omega} \to +\infty$, as $k \to + \infty$. We claim that the second case cannot happen. To this end we will give two proofs of this fact, one works only in dimension $2$, the other one works in any dimension.

\noindent\textbf{Case of dimension 2}: assume that $|\nabla u_{t_k}|_{\infty,\Omega} \to +\infty$. 
Let us set $\Omega^\prime_k:=\{x \in \overline\Omega; \ |\nabla u_{t_k}| \geq \frac{3}{\e}\}$, and $q_k \in \Omega^\prime_k$ such that $|\nabla u_{t_k}|_{\infty,\Omega}=|\nabla u_{t_k}(q_k)|$. We observe that $\Omega_k^\prime$ is closed and hence  is a compact subset of $\overline \Omega$, and being $|\nabla u_t|_{\infty,\partial\Omega} \leq 1-\e$ we have that $\Omega^\prime_k \cap \partial \Omega = \emptyset$ for all $k$. Let $\Omega^{\prime\prime}_k$ be the connected component of $\Omega^\prime_k$ containing $q_k$. Consider now the auxiliary problem

\begin{equation}\label{sistemaauxlap}
\begin{cases}
\Delta v_{t_k} = n t_k (1-\eta_\e^2 |\nabla u_{t_k}|^2) \left(1- \sqrt{1-\eta_\e^2 |\nabla u_{t_k}|^2} e^{ u_{t_k}}  H(e^{ u_{t_k}} q)\right) & \hbox{in} \ \Omega,\\
v_{t_k}=0 & \hbox{in} \ \partial\Omega.
\end{cases}
\end{equation}

We observe that since $u_{t_k}$ is uniformly bounded, then by construction and standard regularity theory we get that $v_{t_k}$ and its gradient are uniformly bounded with respect to $k$. By definition $w_{t_k}:=u_{t_k}-v_{t_k}$ is harmonic in $\Omega^{\prime\prime}_k$. Therefore, considering the isometry $F:\H^2 \to \B^2$, and since harmonicity is preserved through composition with isometries (see \cite[Sect. 2.2]{HeleinWood}), then, $\tilde w_{t_k}:=F \circ w_{t_k}$ is harmonic in $\tilde \Omega^{\prime\prime}_k:= F(\Omega^{\prime\prime}_k) \subset\subset \B^2$. Now, since the hyperbolic metric $\tilde g$ is conformal to the euclidean metric $g_0$ in $\B^2$ (see Sect. 2), we have that $\tilde w_{t_k}$ is harmonic also in $(\tilde \Omega^{\prime\prime}_k,g_0)$. We point out that, in general, this fact is false in the other dimensions. Hence, since $\tilde w_{t_k}$ is harmonic it follows that also $\nabla_0\tilde w_{t_k}$ is harmonic in $\tilde \Omega^{\prime\prime}_k$, so, $\|\nabla_0\tilde w_{t_k}\|_n$ achieves its maximum on the boundary, and thus $\|\nabla_0\tilde w_{t_k}\|_\infty = \|\nabla_0\tilde w_{t_k}\|_{\infty,\partial \tilde \Omega_k^{\prime\prime}}\to +\infty$, as $k\to +\infty$. On the other hand, by construction and \eqref{formulatrasfgrad} we have that
\begin{eqnarray*}
\|\nabla_0\tilde w_{t_k}\|_{\infty, \partial \tilde \Omega^{\prime\prime}_k} &=&\sup_{y \in \partial \tilde \Omega^{\prime\prime}_k} \|\nabla_0\tilde w_{t_k}(y)\|_n=\sup_{q \in \partial \Omega^{\prime\prime}_k}\frac{4}{(1-\|F(q)\|_n^2)^2}|\nabla w_{t_k}(q)|\\
&\leq& \sup_{q \in \partial \Omega^{\prime\prime}_k}\frac{4}{(1-\|F(q)\|_n^2)^2}\left(\frac{3}{\e}+|\nabla v_k(q)|\right)
\end{eqnarray*}
 is uniformly bounded and thus we get a contradiction.

\noindent\textbf{Case of any dimension $n\geq2$}: consider $\tilde u:=u \circ F^{-1}$, where $F:\H^n\to\B^n$ is the hyperbolic stereographic projection. Then $\tilde u$ is a solution of an uniformly elliptic equation which satisfies the hypotheses of \cite[Theorem 15.2]{GT} (see \cite[(i), pag. 367]{GT}). In fact, thanks to Remark \ref{remarkstrutturaoperatoretrasf}, writing $Q^t_\e$ in local coordinates we see by elementary computations that the natural conditions of \cite[(i), pag. 367]{GT}, are satisfied (uniformly in $t$). In particular, introducing the operator $\delta=\frac{\partial}{\partial z} + \sum_{k=1}^n \|p\|^{-2}_n p_k \frac{\partial}{\partial{y_k}}$, we see that $\delta \tilde a_\e^{ij}$, $\delta \tilde b_{\e,t}$ satisfy, as $\|p\|_n\to +\infty$ (uniformly for $(y,z) \in \Lambda \times [\log r_1,\log r_2]$, and in $t \in [0,1]$), the growth conditions of \cite[(15.36)]{GT}, and thus the hypotheses of \cite[Theorem 15.2]{GT} are satisfied with $c\leq0$. 

Thanks to Lemma \ref{lemmaczeroestimates} the oscillation of $u$ is uniformly bounded, moreover, since $|\nabla u|_{\infty,\partial\Omega}\leq 1-\e$ and the structural conditions are satisfied uniformly in $t$, we have that the constant given by  \cite[Theorem 15.2]{GT} is uniformly bounded with respect to $t$. Hence there exists $C$ independent on $t$ such that $\|\nabla_0 \tilde u\|_{\infty,F(\Omega)} \leq C$, and hence, in view of  \eqref{formulatrasfgrad}, the same holds for $|\nabla u|_{\infty,\Omega}$. Therefore, it cannot happen that there exist a sequence $(t_k)$ such that $|\nabla u_{t_k}|_{\infty,\Omega} \to +\infty$, and we are done.

\noindent\textbf{Conclusion}: from the previous discussion the only possibility is that there exists a constant $M$ such that $|\nabla u|_{\infty,\Omega}\leq M$ for all $t \in [0,1]$. From this fact, up to passing to local coordinates, since $Q^t_\e$ is uniformly elliptic (with ellipticity constant independent on $t$), $b_{\e,t}(q,u,\nabla u)$ is uniformly bounded in $t$, thanks to \cite[Theorem 13.7]{GT}, there exists $\alpha_0 \in (0,1)$ and a positive constant $C$, both depending only on $n$, $\Omega$, $|\nabla u|_{\infty,\Omega}$, $\Omega$, and the ratio between the uniform bound on $b_{\e,t}$ and the lower ellipticity constant, such that
$$[\nabla u]_{0,\alpha_0} \leq C,$$
where $[\cdot]_{0,\alpha_0}$ denotes the $C^{0,\alpha_0}$ seminorm. At the end, from this fact and Lemma \ref{lemmaczeroestimates}, we conclude that
$$|u|_{1,\alpha_0} \leq C_1,$$
for some constant $C_1$ not depending on $t$, and the proof is complete.

\end{proof}

\section{Existence and uniqueness of solutions for the regularized problem}
The aim of this section is to prove the following:

 \begin{theorem}\label{teopseudosol}
Let $\alpha \in (0,1)$, $0<r_1 \leq 1 \leq r_2$, with $r_1\neq r_2$, let $\Omega$ be a bounded domain of $\H^n$, with boundary of class $C^{2,\alpha}$. Let $H \in C^{1}({C_{\overline{\Omega}}(r_1,r_2)})$ satisfying the hypotheses $i)$, $ii)$ of Theorem \ref{mainteoprob}.
Assume that $(\Omega,H)$ is admissible. Then, there exists $\bar\e \in (0,1)$ such that for any $\e \in (0,\bar \e)$ Problem \ref{eqpseudoHbumps} has a solution. Moreover such solution is the unique solution of Problem \ref{eqpseudoHbumps} whose associated radial graph is contained in $\mathcal{C}_{\overline\Omega}(r_1,r_2)$. 
\end{theorem}
\begin{proof}
We divide the proof in several steps.\\ 

\noindent\textbf{Step 1}: Choice of $\bar \e \in (0,1)$.\\ 
Let $\hat H$ be the extension of $H$ defined in \eqref{eq:extension}. 
Since $(\Omega,\hat H)$ is admissible we choose $\bar \e=\theta$, where $\theta$ is given by Definition \ref{classadmissible}.\\

 Let $\e \in (0,1)$ such that $\e < \theta$, let $\alpha_0 \in (0,1)$ be the number given by Lemma \ref{lemapriorigradestimates} and set $\beta:=\min\{\alpha,\alpha_0\}$. For any fixed $w  \in C^{1,\beta}(\overline{\Omega})$ we define the operator $L_{w,\e}:C^{2,\beta}_0(\overline\Omega) \to C^{0,\beta}(\overline \Omega)$ as
$$L_{w,\e} u := \sum_{i,j=1}^n \left((1-\eta_\e^2 |\nabla w|^2) \delta_{ij} + \eta_\e^2 w_i w_j\right) u_{ij}.$$

\noindent\textbf{Step 2}: For every  $w \in C^{1,\beta}(\overline{\Omega})$ the operator $L_{w,\e}$ is a bijection of $C^{2,\beta}_0(\overline\Omega)$ onto $C^{0,\beta}(\overline \Omega)$.\\

A mapping $u \in C^{2,\beta}_0(\overline\Omega)$ belongs to the kernel of $L_{w,\e}$ if and only if $u$ solves the Dirichlet problem

\begin{equation}\label{sistemakernelLw}
\begin{cases}
\displaystyle \sum_{i,j=1}^n  \left((1-\eta_\e^2 |\nabla w|^2) \delta_{ij} + \eta_\e^2 w_i w_j\right) u_{ij} = 0 & \hbox{in} \ \Omega,\\[8pt]
u=0 & \hbox{on} \ \partial\Omega.
\end{cases}
\end{equation}

Since $L_{w,\e}$ is uniformly elliptic (see Lemma \ref{lemmaunifelliptic}), by the maximum principle, being $u=0$ on $\partial \Omega$ we obtain that $u=0$ in $\Omega$, and this means that  $L_{w,\e}$ is injective. In order to prove that $L_{w,\e}$ is onto we use the continuity method. Let $t \in [0,1]$, we introduce the family of operators $\mathcal{L}_{t, w, \e}: C^{2,\beta}_0(\overline\Omega) \to C^{0,\beta}(\overline \Omega)$ defined by
$$ \mathcal{L}_{t, w, \e} = (1-t) \Delta + t {L}_{w, \e}.$$
We observe that $\mathcal{L}_{0, w, \e}=\Delta$ and for every $f \in C^{0,\beta}(\overline \Omega)$ the Dirichlet problem

\begin{equation*}
\begin{cases}
\Delta u = f & \hbox{in} \ \Omega,\\
u=0 & \hbox{on} \ \partial\Omega,
\end{cases}
\end{equation*}
admits a solution $C^{2,\beta}(\overline\Omega)$. That is $\mathcal{L}_{0, w, \e}$ sends $C^{0,\beta}(\overline \Omega)$ onto $C^{2,\beta}_0(\overline\Omega)$. Now we claim that there exists a constant $C>0$ such that
\begin{equation}\label{eq1methodcont}
|u|_{2,\beta} \leq C |\mathcal{L}_{t, w, \e}|_{0,\beta},
\end{equation}
for every $t \in [0,1]$, for every $u \in C^{2,\beta}_0(\overline\Omega)$.
In view of the method of continuity this is enough to infer that $\mathcal{L}_{1, w, \e} = L_{w,\e}$ is onto. If \eqref{eq1methodcont} is false then there exist sequences $(t_k) \subset [0,1]$, $(u_k) \subset C^{2,\beta}_0(\overline\Omega)$
such that
\begin{equation}\label{eq2methodcont}
 |\mathcal{L}_{t, w, \e}|_{0,\beta} \to 0 \ \hbox{and} \ |u_k|_{2,\beta} =1.
\end{equation}
By compactness, in particular using also the Ascoli-Arzel\`a Theorem, there exist $t \in [0,1]$ and $u \in C^{2,\beta}_0(\overline\Omega)$ such that, up to a subsequences,

$$ t_k \to t  \ \hbox{and} \ u_k \to u \ \hbox{in} \ C^2(\overline\Omega).$$
By continuity we have $\mathcal{L}_{t, w, \e} u =0$. Since, up to passing in hyperbolic stereographic coordinates, $\mathcal{L}_{t, w, \e}$ is a convex combination of elliptic operators, it is so too. Hence $u=0$. In particular

\begin{equation}\label{eq3methodcont}
u_k \to 0 \ \hbox{in} \ C^0(\overline\Omega).
\end{equation}

We observe that
$$ \mathcal{L}_{t, w, \e} u= \sum_{i,j=1}^n  a^{ij}_{t,\e} u_{ij},$$
where $ a^{ij}_{t,\e}=\left((1-t\eta_\e^2 |\nabla w|^2) \delta_{ij} + t \eta_\e^2 w_i w_j\right)$, and $\mathcal{L}_{t, w, \e}$ is uniformly elliptic, moreover, arguing as in the proof of Lemma \ref{lemmaunifelliptic} we see that the ellipticity constants are independent on $t$. Since the boundary is smooth we can apply Global Schauder estimates and we get that
$$|u_k|_{2,\beta} \leq C (|u_k|_\infty + |\mathcal{L}_{t, w, \e} |_{0,\beta}),$$
with $C$ independent on $k$. This yields a contradiction with \eqref{eq2methodcont}, \eqref{eq3methodcont}. Hence \eqref{eq1methodcont} is true and the proof of Step 2 is complete.\\

\noindent\textbf{Step 3}:  For every $C>0$ there exists $K>0$ such that if $|w|_{1,\beta} \leq C$ then $|u|_{2,\beta} \leq K |L_{w,\e}u|_{0,\beta}$ for every $u \in C^{2,\beta}_0(\overline \Omega)$.\\
We argue by contradiction as in the last part of the proof of Step 2. If the result is false then there exist a bounded sequence $(w_k)$ in $C^{1,\beta}(\overline \Omega)$ and a sequence $(u_k)$ in $C^{2,\beta}_0(\overline \Omega)$ such that
\begin{equation}\label{eq4methodcont}
|u_k|_{2,\beta}=1 \ \hbox{and} \ |L_{w_k,\e} u_k|_{0,\beta} \to 0.
\end{equation}
By compactness, there exist $w \in C^1(\overline \Omega)$ and $u \in C^2_0(\overline\Omega)$ such that, up to subsequences,
$$ w_k \to w \ \hbox{in} \ C^1(\overline \Omega) \ \hbox{and} \ u_k \to u \ \hbox{in} \ C^2(\overline\Omega).$$
By continuity we get $L_{w,\e} u=0$. Then $u=0$, by Step 2. Taking into account of Lemma \ref{lemmaunifelliptic} we observe that the operators $L_{w_k, \e}$ are uniformly elliptic with ellipticity constants independent on $k$. Using standard Schauder's estimates we obtain that
$$ |u_k|_{2,\beta} \leq C_1 \left( |u_k|_\infty + |L_{w_k, \e}u_k|_{0,\beta}\right),$$
where $C_1$ is constant independent on $k$. Since $u_k \to 0$ in $C^0(\overline \Omega)$ and by \eqref{eq4methodcont} we reach a contradiction. The proof of Step 3 is complete.\\

\noindent\textbf{Step 4}: Let $(w_k)$ be a bounded sequence in $C^{1,\beta}(\overline \Omega)$ and let $(f_k)$ be a bounded sequence in  $C^{0, \beta}(\overline \Omega)$. Then the sequence $(u_k)$ of solutions of
\begin{equation}\label{sistmethodcont5}
\begin{cases}
L_{w_k, \e} u_k = f_k & \hbox{in} \ \Omega,\\
u_k=0 & \hbox{on} \ \partial\Omega.
\end{cases}
\end{equation}
is bounded in $C^{2,\beta}(\overline \Omega)$.\\[8pt]
The existence of a solution $u_k$ of \eqref{sistmethodcont5} is given by Step 2, and the thesis follows from Step 3.\\


\noindent\textbf{Step 5}: 
Let us consider the map $T_\e: C_0^{1,\beta}(\overline \Omega) \to C^{1,\beta}_0(\overline{\Omega})$, defined as follows: for every $w \in C_0^{1,\beta}(\overline \Omega)$ we set $T_{\e} (w):= u$, where $u=u(w,\e)$ is the unique solution of the problem

\begin{equation}\label{sistemaaux1}
\begin{cases}
\displaystyle \sum_{i,j=1}^n \left((1-\eta_\e^2 |\nabla w|^2) \delta_{ij} + \eta_\e^2 w_i w_j\right) u_{ij} = n  (1-\eta_\e^2 |\nabla w|^2) \left(1- \sqrt{1-\eta_\e^2 |\nabla w|^2} e^w  \hat H(e^w q)\right) & \hbox{in} \ \Omega,\\[6pt]
u=0 & \hbox{on} \ \partial\Omega.
\end{cases}
\end{equation}
We claim that $T_\e$ is a compact operator.\\[8pt]
We first observe $T_{\e}$ is well defined, in fact, as proved in Step 2, for a given $w \in C_0^{1,\beta}(\overline \Omega)$, the operator $L_{w,\e}$ is a bijection between $C^{2,\beta}_0(\overline \Omega)$ and $C^{0,\beta}(\overline \Omega)$. In addition is $T_{\e}$ is a linear map. It remains to prove that $T_\e$ maps bounded families of $C_0^{1,\beta}(\overline \Omega)$ into relatively compact subsets of $C_0^{1,\beta}(\overline \Omega)$.

Let $(w_\lambda)$ be a bounded family of $C_0^{1,\beta}(\overline \Omega)$, then $(u_\lambda)$, where $u_\lambda=Tw_\lambda$, is a family of solutions of \ref{sistemaaux1}. Hence $u_\lambda \in C^{2,\beta}_0(\overline\Omega)$ and since we assuming that there exists $C>0$ such that $|w_\lambda|_{1,\beta}\leq C$ then by Step 3 we have
$$|u_\lambda|_{2,\beta} \leq K \left|n (1-\eta_\e^2 |\nabla w_\lambda|^2) \left(1- \sqrt{1-\eta_\e^2 |\nabla w_\lambda|^2} e^{w_\lambda}  \hat H(e^{w_\lambda} q)\right)\right|_{0,\beta} \leq K_1,$$
where $K_1$ is a positive constant not depending on the family. Hence $(u_\lambda)$ is uniformly bounded in $C^{2,\beta}_0(\overline\Omega)$, and in particular by Step 4 and Ascoli-Arzel\`a Theorem it is relatively compact in $C^{1,\beta}_0(\overline\Omega)$. This proves that $(u_\lambda)$ is relatively compact in $C^{1,\beta}_0(\overline\Omega)$ and we are done.\\

\noindent\textbf{Step 6}: There exists a constant $C>0$ such that $|u|_{1,\beta} \leq C$ for any $u \in C^{1,\beta}_0(\overline\Omega)$ satisfying $u=tT_\e u$, where $t \in [0,1]$.\\[8pt]
We first observe that by definition and standard elliptic regularity theory any  $u \in C^{1,\beta}_0(\overline\Omega)$ satisfying $u=tT_\e u$ is of class $C_0^{2,\beta}(\overline\Omega)$ and verifies $\hat Q^t_\e(u)=0$. 
Thanks to Lemma \ref{lemapriorigradestimates}, and since $\beta \leq \alpha_0$, there exists $C>0$ such that $|u|_{1,\beta}\leq C$, provided that $|\nabla u|_{\infty,\partial\Omega}\leq1-\e$. Therefore, in oder to conclude it is sufficient to check this boundary estimate for the gradient.

Let $q_0 \in \partial\Omega$ such that $|\nabla u(q_0)|=|\nabla u|_{\infty,\partial\Omega}$. If $\nabla u(q_0)=0$ it follows that $\nabla u=0$ on $\partial \Omega$ and hence there is nothing to prove. Therefore, let us assume that $\nabla u(q_0)\neq 0$. 

Since $(\Omega,\hat H)$ is admissible, for any $t \in [0,1]$ there exist $\varphi_1, \varphi_2 \in C^{2}(\overline\Omega)$ satisfying i)-iv) of Definition \ref{classadmissible} at $q_0$. Hence, taking into account of the choice of $\e$ and Remark \ref{remoperator}, we have
$$  \hat Q^t_{\e}(\varphi_1) \geq   \hat Q^t_\e (u) \geq  \hat Q^t_\e(\varphi_2) \ \ \hbox{in} \ \Omega,$$
and $\varphi_1 \leq u \leq \varphi_2$ on $\partial \Omega$.
Let us write $\hat Q_\e^t(u)=\sum_{i,j=1}^n a_\e^{ij} u_{ij} +   \hat b_{\e,t}(q,u,\nabla u)$, where $$ \hat b_{\e,t}(q,u,\nabla u):=- n t (1-\eta_\e^2 |\nabla u|^2)  + n t(1-\eta_\e^2 |\nabla u|^2)^{3/2} e^u  \hat H(e^u q).$$ 
Notice that thanks to assunption ii) and Remark \ref{remextensadmissible} it follows that for any fixed $q \in \Omega$ the map $z \mapsto e^z  \hat H(e^zq)$ is not increasing. 

Thanks to  Lemma \ref{lemmaunifelliptic} and  Remark \ref{remarkstrutturaoperatoretrasf}, under the hyperbolic stereographic projection $F:\H^n \to \B^n$, the operator $\hat Q_\e^t$ is transformed into an uniformly elliptic operator of the form $$\tilde Q^t_\e \tilde u= \sum_{i,j=1}^n \tilde a_\e^{ij}(y,\nabla_0 \tilde u ) \tilde u_{ij} +  \tilde b_{\e,t}(y,\tilde u,\nabla_0 \tilde u),$$ where $y \in F(\Omega)$, $\nabla_0 \tilde u$ is the euclidean gradient, $\tilde u_{ij}$ are the second partial derivatives of $\tilde u=u\circ F^{-1}$. In view of Remark \ref{remarkstrutturaoperatoretrasf} and assumption ii) the principal part $\tilde a_\e^{ij}(y,p)$ does not depend on $z$, and for each $(y,p) \in F(\Omega) \times \R^n$ the map $z \mapsto  \tilde b_{\e,t}(y,z,p)$ is non-increasing. Hence the comparison principle applies (see \cite[Theorem 10.1]{GT}), and thus setting $\tilde \varphi_i:=\varphi_i \circ F^{-1}$, for $i=1,2$, from $\tilde Q^t_\e(\tilde \varphi_1) \geq \tilde Q^t_\e(\tilde u) \geq \tilde Q^t_\e(\tilde \varphi_2)$ in $F(\Omega)$, and $\tilde\varphi_1 \leq \tilde u  \leq \tilde\varphi_2$ on $\partial F(\Omega)$, it follows that $\tilde\varphi_1 \leq \tilde u  \leq \tilde\varphi_2$ in $F(\Omega)$. Therefore we obtain
\begin{equation}\label{stimau}
\varphi_1 \leq u \leq \varphi_2\ \hbox{in}\ \Omega.
\end{equation}

We observe that since $u=0$ on $\partial\Omega$ then $\nabla u(q_0)$ is orthogonal to $T_{q_0}\partial\Omega$, where $T_{q_0}\partial\Omega$ is the tangent space at $q_0$ for $\partial\Omega$, and we have the orthogonal decomposition $Span\{\nabla u(q_0)\} \bigoplus T_{q_0} \partial\Omega =  T_{q_0}\H^n$.

Let us set $ \hat w:=\frac{\nabla u(q_0)}{|\nabla u(q_0)|}$ and consider a curve  $\psi:(-\delta,\delta) \to \H^n$ such that $\psi(0)=q_0$,  $\psi(s) \in \Omega$ for $s \in (0,\delta)$ and $\psi^\prime(0)=\hat w$ if $\hat w$ points towards the interior of $\Omega$ (otherwise we take $\psi^\prime(0)=- \hat w$). Since $\Omega$ has a smooth boundary we can always find a curve satisfying these properties. From \eqref{stimau}, and since $u(q_0)= \varphi_1(q_0)= \varphi_2(q_0)=0$ we deduce that for all sufficiently small $h>0$

\begin{equation}\label{stimau2}
\frac{\varphi_1 (\psi(h)) - \varphi_1(\psi(0))}{h} \leq \frac{u(\psi(h)) - u(\psi(0))}{h} \leq \frac{\varphi_2(\psi(h)) - \varphi_2(\psi(0))}{h}.
\end{equation}
Passing to the limit as $h\to 0^+$ we get that

$$ d \varphi_1(q_0)[\hat w] \leq du(q_0)[\hat w] \leq d \varphi_2(q_0)[\hat w],$$
(if $\hat w$ points in the opposite direction \eqref{stimau2} holds but with the reversed inequalities). At the end, it follows that
$$| du(q_0)[\hat w]| \leq \max\{|d \varphi_1(q_0)[\hat w]|, |d \varphi_2(q_0)[\hat w]|\}.$$
Since $\H^n$ a spacelike hypersurface, then for any $q \in \H^n$ the Cauchy-Schwartz inequality holds in $T_q\H^n$ for $\l\cdot,\cdot\r_{Tq\H^n}$ (we point out that, in general, the Cauchy-Schwartz inequality does not hold  in $\L^{n+1}$, see \cite{Lopez2013}). In particular $|d \varphi_i(q_0)[\hat w]|=|\l \nabla \varphi_i(q_0), \hat w \r| \leq | \varphi_i(q_0)| | \hat w| = | \varphi_i(q_0)|$.

Hence, by the previous discussion and by Definition \ref{classadmissible} we have
$$|\nabla u(q_0)|=|\l\nabla u (q_0), \hat w\r| = |du(q_0)[\hat w]| \leq \max\{|d \varphi_1 [\hat w], d \varphi_2 [\hat w] \}\leq  1-\theta.$$ 

At the end since $\e < \theta$ we get that $|\nabla u|$ verifies the desired boundary estimate, and thus from the initial discussion, the proof of Step 6 is complete.\\

\noindent\textbf{Step 7:} Existence of a solution of Problem \eqref{eqpseudoHbumps}.\\[8pt]
Thanks to Step 5 and Step 6 it follows that the operator $T_\e: C^{1,\beta}_0(\overline\Omega) \to C^{1,\beta}_0(\overline\Omega)$ satisfies the hypotheses of the Leray-Schauder Theorem (see \cite[Theorem 11.3]{GT}), and thus there exists $u \in C^{1,\beta}_0(\overline\Omega)$ which solves $u=T_\e u$. Hence, $u \in C^{2,\beta}_0(\overline\Omega)$, and by definition of $\hat H$ and Lemma \ref{lemmaczeroestimates} $u$ is a solution of  Problem \eqref{eqpseudoHbumps}. The proof of Step 7 is complete.\\

\noindent\textbf{Step 8:} Uniqueness.\\[8pt]
For the uniqueness of the solution it is sufficient to argue as in \cite[Sect. 2.3]{CaldGull}. For the sake of completeness we give a sketch of the proof.

Let us fix $\e \in (0,1)$ and let $u_1,u_2\in C^{2}(\Omega)\cap C^{0}(\overline\Omega)$ be two solutions of Problem \ref{eqpseudoHbumps} such that the corresponding radial graphs are contained in $\mathcal{C}_{\overline\Omega}(r_1,r_2)$. 
If $u_1\neq u_2$ then there exists $\overline{q}\in\Omega$ such that $u_1(\overline{q})\ne u_2(\overline{q})$. Without loss of generality we can assume that $u_1(\overline{q})<u_2(\overline{q})$. Then there exists $\mu>0$ such that $u_1(q)+\mu\ge u_2(q)$ for every $q\in\Omega$ and $u_1(q_{0})+\mu=u_2(q_{0})$ at some $q_{0}\in\Omega$. Set $\overline{u_1}:=u_1+\mu$ and observe that $\overline{u_1}$ verifies
$$
\sum_{i,j=1}^n\left((1-\eta_\e^2|\nabla\overline{u_1}|^{2})\delta_{ij}+\eta_\e^2\overline{u_1}_{i}\overline{u_1}_{j}\right)\overline{u_1}_{ij}\le n(1-|\eta_\e^2\nabla \overline{u_1}|^{2})\left(1-\sqrt{1-\eta_\e^2|\nabla \overline{u_1}|^{2}}~\!e^{\overline{u_1}}~\!\hat H(e^{\overline{u_1}}q)\right)\quad\textrm{in $\Omega$}
$$
because of ii) and $\mu>0$. Notice that the radial graph defined by $\overline{u}_1$ stays over (in the radial direction) that one corresponding to $u_2$ and they touch themselves at the point $X_{0}=q_{0}e^{u_2(q_{0})}$. Now, in order to conclude it sufficient to compare $\overline{u_1}$ and $u_2$ by means of the Hopf maximum principle. To this end we use the version stated in \cite[Theorem 2.3]{PuSe04} for the operator 
\begin{equation*}
\begin{split}
 Q_\e(u)=&\sum_{i,j=1}^n \left((1-\eta_\e^2|\nabla u|^{2})\delta_{ij}+u_{i}u_{j}\right)u_{ij}\\
&\qquad-n(1-\eta_\e^2|\nabla u|^{2})\left(1-\sqrt{1-\eta_\e^2|\nabla u|^{2}}~\!e^{u}~\! \hat H(e^{u}q)\right).
\end{split}
\end{equation*}
It is easy to see that, up to passing to hyperbolic stereographic coordinates, the assumptions of \cite[Theorem 2.3]{PuSe04} are fulfilled, and applying the theorem as in \cite{CaldGull}, we deduce that $\overline{u}_1=u_2$ in $\Omega$. But this gives a contradiction since $\overline{u_1}|_{\partial\Omega}=\mu>0=u_2|_{\partial\Omega}$. Hence it must be $u_1=u_2$ and we are done.\\


\end{proof}

\section{An interior estimate for the gradient}

In this section we prove an estimate for the gradient when the maximum point of its modulus lies in the interior of the domain $\Omega$. We begin with a preliminary elementary result of linear algebra. 

\begin{lemma}\label{lemmaalgebralinear}
Let $A=(a_{ij}),B=(b_{ij}) \in \mathcal{M}_n(\R)$ be two symmetric matrices. Assume that $A$ is positive semi-definite and $B$ is negative semi-definite. Then $$\sum_{i,j=1}^n a_{ij} b_{ij} \leq 0.$$
\end{lemma}

\begin{proof}
Since $A$, $B$ are symmetric we have that $\sum_{i,j=1}^n a_{ij} b_{ij}=trace(AB)$, and there exist two invertible matrices $P$, $Q$ such that $P^{-1}AP=D_A$, $Q^{-1}BQ=D_B$ are diagonal. Thanks to the assumptions we have that $D_A$ has non-negative elements on the diagonal, while $D_B$ has non positive elements on the diagonal. Therefore, since the trace is invariant under similitude, and diagonal matrices commute in the product, we have
$$trace(AB)=trace(P^{-1}APP^{-1}BP)=trace(D_AP^{-1}BP)=trace(P^{-1}D_ABP)=trace(D_AB).$$
Now, by the same argument we get that
$$trace(D_AB)=trace(Q^{-1}D_AQQ^{-1}BQ)=trace(D_AQ^{-1}QD_B)=trace(D_AD_B).$$
Therefore, $\sum_{i,j=1}^n a_{ij} b_{ij}=trace(AB)=trace(D_AD_B)\leq 0$, and the proof is complete.
\end{proof}

\begin{definition}\label{tangentialcomponent}
Let  $H \in C^{1}(\mathcal{C}_{\overline\Omega})$, let $\nabla_0 H$ be gradient of $H$ in $\R^{n+1}$ with respect to the flat metric. We define the (euclidean) tangential component of $\nabla_0 H$ on $T_{\frac{x}{|x|}} \H^n$ as the vector
$$\nabla^T_0 H (x):=\nabla_0 H(x) - (\nabla_0 H(x), \hat r(x))_{n+1}\ \hat r(x), \ \ x \in \mathcal{C}_{\overline\Omega}$$
where $$\hat r(x):=\frac{(x_1,\ldots,x_n,-x_{n+1})}{\|x\|_{n+1}}.$$ 
\end{definition}

\begin{remark}\label{rmk:tangentialcomponent}
We point out that by definition $\nabla^T_0 H (x)=\nabla_0 H(x)- \l \nabla_0H(x), \frac{x}{\|x\|_{n+1}}\r \hat r(x)$, and if $v \in \R^{n+1}$ is such that $(\hat r(x), v)_{n+1}=0$ then $\l \frac{x}{|x|}, v\r=0$, and viceversa. In particular $(\hat r(x), v)_{n+1}=0$ for any $v \in T_{\frac{x}{|x|}} \H^n$, $x \in \mathcal{C}_{\overline\Omega}$
\end{remark}

In the sequel we will make use also of the following formulas for the second and third covariant derivatives of a smooth function $u$ defined over $\Omega$ (see \cite[Sect. 2]{Yau})
\begin{equation}\label{formulasecondcovder}
\sum_{j=1}^nu_{ij} \omega^j= du_i - \sum_{j=1}^nu_j \omega_{ji},
\end{equation}

\begin{equation}\label{formulathirdcovder}
\sum_{k=1}^nu_{ijk} \omega^k = du_{ij} - \sum_{k=1}^n u_{kj} \omega_{ki} - \sum_{k=1}^nu_{ik} \omega_{kj}.
\end{equation}
\begin{proposition}\label{stimagradint}
Let $\Omega$ be a bounded domain of $\H^n$, let $H \in C^{1}(\mathcal{C}_{\overline\Omega})$, $\e \in (0,1)$ and let $u \in C^{3}(\overline{\Omega})$ be a solution of
\begin{equation}\label{eq:revfinal}  
\sum_{i,j=1}^n((1- \eta^2_\e|\nabla u|^2)\delta_{ij}+\eta^2_\e u_iu_j) u_{ij}= n (1-\eta_\e^2 |\nabla u|^2) \left(1- \sqrt{1-\eta_\e^2 |\nabla u|^2} e^u H(e^u q)\right).
\end{equation}
 Then, if the maximum point $q_0$ of $|\nabla u|$ lies in the interior of $\Omega$, we have $|\nabla u(q_0)|=0$ or
\begin{equation}\label{eqlemmagradient}
\begin{array}{lll}
&&\displaystyle \left[- (n-1) -n (1-\eta^2_\e |\nabla u(q_0)|^2)^{1/2} e^{u(q_0)}  \frac{\partial}{\partial \lambda} \left(\lambda H(\lambda q)\right){\big|_{\lambda=e^{u(q_0)}}}\right] |\nabla u(q_0)|\\[6pt]
&&\displaystyle \ \ \ \ \ \ \ \ \ \ \ \ \   - n^{3/2} (1-\eta^2_\e |\nabla u(q_0)|^2)^{1/2} e^{2u(q_0)} \|\nabla_0^T H(e^{u(q_0)}q_0)\|_{n+1} \ \leq 0,
\end{array}
\end{equation}
where $\nabla_0^T H$ is the (euclidean) tangential component of $(\nabla_0 H)(e^{u(q_0)}q_0)$ on $T_{q_0} \H^n$.
\end{proposition}

\begin{proof}
We will prove a more general version of \eqref{eqlemmagradient}. 
Let us fix a smooth positive function $f:\R \to \R^+$ and consider the auxiliary function $\varphi:=f(2Cu)|\nabla u|^2$, where $C \in \R$ is a fixed constant. In order to simplify the notation we set $v:=|\nabla u|^2$, hence $\varphi=f(2Cu) v$.  Assume that $\varphi$ has a maximum point at some $q_0$ lying in the interior of $\Omega$. Hence $\nabla \varphi (q_0)=0$ and the Hessian $(\varphi_{ij}(q_0))$ is negative semi-definite.

By direct computation we have $v_i = 2\sum_{j=1}^n u_j u_{ji}$ and from $\nabla \varphi (q_0)=0$ we get that

\begin{equation}\label{eqT15}
\sum_{h=1}^nf(2Cu)u_hu_{hi} + f^\prime(2Cu)Cvu_i=0, \ \hbox{for all} \ i=1,\ldots,n,
\end{equation}
which implies that

\begin{equation}\label{eqT16}
\sum_{i,h=1}^n u_iu_{ih}u_h =- C \frac{f^\prime}{f}v^2,
\end{equation}
where, $f$, $f^\prime$ stand, respectively, for $f(2Cu)$, $f^\prime(2Cu)$. By a simple computation, from \eqref{eqT16}, we get that
\begin{equation}\label{eqT17}
\sum_{i,h,k=1}^n u_iu_{ih}u_{hk}u_k = \left(\frac{f^\prime}{f}\right)^2 C^2v^3.
\end{equation}

Let us set $ a_\e^{ij}:=((1- \eta^2_\e|\nabla u|^2)\delta_{ij}+\eta^2_\e u_iu_j)$, $b_\e:=n (1-\eta_\e^2 |\nabla u|^2) \left(1- \sqrt{1-\eta_\e^2 |\nabla u|^2} e^u H(e^u q)\right)$. 
Since $(a_\e^{ij})$ is a positive definite symmetric matrix and $(\varphi_{ij}(q_0))$ is symmetric negative semi-definite, then, from Lemma \ref{lemmaalgebralinear} it follows that
\begin{equation}\label{eqT14bis}
 \sum_{i,j=1}^n a_\e^{ij}\varphi_{ij}(q_0) \leq 0.
\end{equation}

In order to get an estimate for $v=|\nabla u|^2$ the idea is to use \eqref{eqT14bis}. To this end we compute explicitly $\varphi_{ij}(q_0)$.
Recalling that $\varphi_i= 2  \left(\sum_{h=1}^nf(2Cu)u_hu_{hi} + f^\prime(2Cu)Cvu_i \right)$, and using \eqref{eqT15}, \eqref{formulasecondcovder} we have

\begin{eqnarray*}
&&\sum_{j=1}^n \varphi_{ij}(q_0) \omega^j \\
&=& 2  \sum_{j=1}^n \left[ \sum_{h=1}^n \left(2C  f^\prime u_j u_h u_{hi} + f u_{hj} u_{hi} +  f u_{hij} u_h \right)+ 2C^2 f^{\prime\prime} u_j u_i v + Cf^\prime u_{ij} v +  \sum_{h=1}^n 2C f^\prime u_i u_h u_{hj}\right] \omega^j,
\end{eqnarray*}
and thus from \eqref{eqT14bis} we infer that
\begin{equation}\label{eqT14}
\begin{array}{lll}
 &&\displaystyle 2 \left[ \sum_{i,h=1}^n \left(4C f^\prime (1-\eta^2_\e v) u_i u_h u_{hi} + f(1-\eta^2_\e v)u_{hi}^2 + f (1-\eta^2_\e v)u_{hii}u_h \right) +  2C^2 f^{\prime\prime} (1-\eta^2_\e v) v^2 \right.\\[10pt] 
 &&\ \displaystyle +  \sum_{i=1}^n C f^\prime (1-\eta^2_\e v) u_{ii} v +   \sum_{i,h=1}^n 2C f^\prime  \eta^2_\e u_i u_h u_{hi} +  \sum_{i,j,h=1}^n \left( f \eta^2_\e  u_i u_j u_{hj} u_{hi} + f \eta^2_\e u_i u_j u_h u_{hij}\right) \\[4pt]
 &&\ \displaystyle \left. +  2C^2 f^{\prime\prime} \eta v^3+  \sum_{i,j=1}^n C f^\prime \eta^2_\e u_i u_j u_{ij} v +  \sum_{j,h=1}^n 2C f^\prime \eta^2_\e  u_j  u_ h u_{hj}  v \right] \leq 0.
\end{array}
\end{equation}

Now we estimate and rewrite the terms involving the second and third covariant derivatives.
We may choose a coordinate frame at $q_0$ satisfying $\delta_{1i}v^{1/2}=u_i$. Either $v(q_0)=0$ in which case $\max_{\overline{\Omega}}\varphi = \varphi(q_0)=0$ and the thesis follows immediately. Otherwise in these coordinates, from \eqref{eqT15}, it follows that
\begin{equation}\label{eqT18}
u_{11} =- \frac{f^\prime}{f}Cv,
\end{equation}
which implies that
\begin{equation}\label{eqT19}
\sum_{i,j=1}^n u_{ij} u_{ij} \geq  \left(\frac{f^\prime}{f}\right)^2 C^2v^2.
\end{equation}
Since $u$ is a solution of \eqref{eq:revfinal}, computing at $q_0$ in these coordinates we have
$$\sum_{i=1}^n(1-\eta^2_\e v) u_{ii} + \eta^2_\e v u_{11} = b_\e,$$
and from \eqref{eqT18} we obtain 
\begin{equation}\label{eqT20}
\sum_{i=1}^n u_{ii} = \left( b_\e +C \frac{f^\prime}{f}\eta^2_\e v^2 \right)  (1-\eta^2_\e v)^{-1}.
\end{equation}

Recalling that $\eta_\e(|\nabla u|)=\eta_\e(v^{1/2})$, by direct computation we infer
\begin{equation}\label{eq:revetanabla}
\nabla_k(\eta^2_\e v)=  \sum_{h=1}^n \left(2 \eta_\e \eta^\prime_\e v^{1/2}u_h u_{hk} + 2 \eta^2_\e u_h u_{hk}\right),
\end{equation}
where it is understood that $\eta^\prime_\e$ stands for $\eta^\prime_\e(v^{1/2})$.
By differentiating \eqref{eq:revfinal}, taking into account of \eqref{formulasecondcovder},\eqref{formulathirdcovder} and \eqref{eq:revetanabla}, after some standard computations we deduce that

\begin{equation*}
\begin{array}{lll}
&&\displaystyle \sum_{k=1}^n \left[ - \sum_{i,h=1}^n (2 \eta_\e \eta^\prime_\e v^{1/2}u_h u_{hk} + 2 \eta^2_\e u_h u_{hk}) u_{ii}+ \sum_{i,j,h=1}^n2 \eta_\e \eta^\prime_\e v^{-1/2}u_i u_j u_h u_{hk} u_{ij} \right.\\[4pt]
 && \  \left.\displaystyle + \sum_{i,j=1}^n \left(\eta^2_\e u_j u_{ik} u_{ij} + \eta^2_\e u_i u_{jk} u_{ij}  + ((1-\eta^2_\e v)\delta_{ij}+\eta^2_\e u_i u_j) u_{ijk}\right)\right] \omega^k =\displaystyle \sum_{k=1}^n  (b_\e)_k\omega^k.
\end{array}
\end{equation*}
Now, contracting the equation with $u_k$, we get that

\begin{equation}\label{eqT21}
\begin{array}{lll}
&& \displaystyle \sum_{i,h,k=1}^n \left(- 2 \eta_\e \eta^\prime_\e v^{1/2}u_h u_k u_{hk} u_{ii} - 2 \eta^2_\e u_h u_k u_{hk} u_{ii} \right)+ \sum_{i,j,h,k=1}^n 2 \eta_\e \eta^\prime_\e v^{-1/2}u_i u_j u_k u_h u_{hk} u_{ij} \\[4pt]
 && \ \displaystyle +  \sum_{i,j,k=1}^n \left( \eta^2_\e u_j u_ k u_{ik} u_{ij} +  \eta^2_\e u_i u_k u_{jk} u_{ij} \right) +  \sum_{i,k=1}^n(1-\eta^2_\e v) u_k u_{iik}+  \sum_{i,j,k=1}^n \eta^2_\e u_i u_j u_k u_{ijk}\\
 & =& \displaystyle \sum_{k=1}^n  (b_\e)_k u_k.
\end{array}
\end{equation}

Since the Ricci curvature of the Hyperbolic space is $R_{ij}=-(n-1)\delta_{ij}$, the Ricci formula (see formula (2.11) in \cite{Yau}) gives

\begin{equation}\label{eqT22}
\sum_{k=1}^n u_k u_{kii}= \sum_{k=1}^n u_k u_{iik} - (n-1) v.
 \end{equation}
 
Hence, using \eqref{eqT22}, and taking into account of \eqref{eqT16}, \eqref{eqT17},\eqref{eqT20}, we rewrite \eqref{eqT21} as

\begin{equation}\label{eqT23}
\begin{array}{lll}
&& \displaystyle \sum_{i,k=1}^n (1-\eta^2_\e v) u_k u_{kii}+ \sum_{i,j,k=1}^n \eta^2_\e u_i u_j u_k u_{ijk}\\[12pt]
&=& - (n-1) v (1-\eta^2_\e v)  - 2 \eta_\e \eta^\prime_\e C  \frac{f^\prime}{f} v^{5/2} (1-\eta^2_\e v)^{-1} (b_\e+C \frac{f^\prime}{f}\eta^2_\e v^2)  \\[8pt]
 && - 2 \eta^2_\e  C v^{2} (1-\eta^2_\e v)^{-1} (b_\e+C\eta^2_\e v^2)\\[4pt]
  &&\displaystyle - 2\eta_\e \eta^\prime_\e C^2 \left( \frac{f^\prime}{f}\right)^2 v^{7/2} - 2 \eta^2_\e C^2 \left( \frac{f^\prime}{f}\right)^2 v^3 + \sum_{k=1}^n (b_\e)_k u_k.
\end{array}
\end{equation}

Now, from \eqref{eqT14} and by using \eqref{eqT16}, \eqref{eqT17},  \eqref{eqT19}, \eqref{eqT20} \eqref{eqT23}, we deduce that

\begin{equation}\label{eqT24}
\begin{array}{lll}
&& - 4C^2 (1-\eta^2_\e v) \frac{(f^\prime)^2}{f} v^2+ C^2 (1-\eta^2_\e v) \frac{(f^\prime)^2}{f} v^2 \\[6pt]
&&+ f\left[ - (n-1) v (1-\eta^2_\e v)  - 2 \eta_\e \eta^\prime_\e C  \frac{f^\prime}{f} v^{5/2} (1-\eta^2_\e v)^{-1} (b_\e+C  \frac{f^\prime}{f}\eta^2_\e v^2)  \right.\\[6pt]
&&\displaystyle  \left. - 2 C \eta^2_\e   v^{2} (1-\eta^2_\e v)^{-1} (b_\e+C \frac{f^\prime}{f}\eta^2_\e v^2) - 2 C^2\eta_\e \eta^\prime_\e   \left(\frac{f^\prime}{f}\right)^2 v^{7/2}- 2 C^2  \left(\frac{f^\prime}{f}\right)^2 v^{3} + \sum_{k=1}^n (b_\e)_k u_k\right] \\[10pt]
&& + 2C^2 f^{\prime\prime} (1-\eta^2_\e v) v^2 + C f^\prime v (b_\e+C  \frac{f^\prime}{f}\eta^2_\e v^2) - 2C^2 \frac{(f^\prime)^2}{f} \eta^2_\e v^2  + C^2 \frac{(f^\prime)^2}{f} \eta^2_\e v^3 \\[8pt]
&&+ 2C^2f^{\prime\prime} \eta^2_\e v^3 - C^2 \frac{(f^\prime)^2}{f} \eta^2_\e v^3 - 2C^2 \frac{(f^\prime)^2}{f} \eta^2_\e v^3\leq 0.
\end{array}
\end{equation}

Now we compute the term $ \sum_{k=1}^n (b_\e)_k u_k$. To this end let us observe that $$\left(e^u H(e^uq)\right)_k(q_0) = e^{u(q_0)} u_k H(e^{u(q_0)}q_0)+e^{2{u(q_0)}}u_k(q_0)(\nabla_0H)(e^{u(q_0)}q_0) \cdot q_0 + e^{2{u(q_0)}}(\nabla_0H)(e^{u(q_0)}q_0) \cdot e_k(q_0),$$ where $\cdot$ denotes the standard euclidean product of $\R^{n+1}$, $\nabla_{0} H$ is the gradient of $H$ with respect to the flat metric in $\R^{n+1}$. Then, after some computations and taking into account of \eqref{eq:revetanabla}, we get that 
\begin{equation}\label{eqT24new}
\begin{array}{lll}
\displaystyle \sum_{k=1}^n (b_\e)_k u_k(q_0)&=& \displaystyle 2nC  \frac{f^\prime}{f}  \eta_\e \eta^\prime_\e  v^{5/2} \left(1-\sqrt{1-\eta^2_\e v} \ e^{u(q_0)} H(e^{u(q_0)}q_0)\right) \\[8pt]
&+&  \displaystyle   2nC    \frac{f^\prime}{f} (\eta^2_\e v)   v \left(1-\sqrt{1-\eta^2_\e v} \ e^{u(q_0)} H(e^{u(q_0)}q_0)\right) \\[8pt]
&-&  \displaystyle   n C   \frac{f^\prime}{f}  \eta_\e \eta^\prime_\e \sqrt{1-\eta^2_\e v} \left(e^{u(q_0)} H(e^{u(q_0)}q_0)\right) v^2 \\[8pt]
&-&  \displaystyle   n C   \frac{f^\prime}{f}  (\eta^2_\e v) (1-\eta^2_\e v) \left(e^{u(q_0)} H(e^{u(q_0)}q_0)\right) v \\[8pt]
&-&  \displaystyle   n    (1-\eta^2_\e v)^{3/2} e^{u(q_0)}  \left(e^{u(q_0)} H(e^{u(q_0)} q_0)+ e^{2u(q_0)} \nabla_{0} H(e^{u(q_0)} q_0) \cdot q_0\right) v \\[8pt]
&-&  \displaystyle   n   (1-\eta^2_\e v)^{3/2} e^{2u(q_0)}  \nabla_{0} H(e^{u(q_0)}q_0) \cdot \nabla u,
\end{array}
\end{equation}
Computing $ \frac{\partial}{\partial \lambda} \left(\lambda H(\lambda q_0)\right){\big|_{\lambda=e^{u(q_0)}}}$ and taking into account of Remark \ref{rmk:tangentialcomponent}, we rewrite the last two terms of \eqref{eqT24new} as
\begin{equation}\label{eqT24bis}
-  \displaystyle   n    (1-\eta^2_\e v)^{3/2} e^{u(q_0)} \frac{\partial}{\partial \lambda} \left(\lambda H(\lambda q_0)\right){\big|_{\lambda=e^{u(q_0)}}} v 
-  \displaystyle   n   (1-\eta^2_\e v)^{3/2} e^{2u(q_0)}  \nabla_0^T H(e^{u(q_0)}q_0) \cdot \nabla u,
\end{equation}
where $\nabla_0^T H$ is the euclidean tangential component of $(\nabla_0 H)$ on $T_{q_0} \H^n$.
Hence, from  \eqref{eqT24}, \eqref{eqT24new}, \eqref{eqT24bis}, regrouping and simplifying terms, we get that

\begin{equation}\label{eqT25}
\begin{array}{lll}
&& \displaystyle - 2C^2 (f^\prime)^2 \left[1-\eta^2_\e + \eta^\prime_\e (\eta_\e v^{1/2})\right] v^3\\[8pt] 
 && \displaystyle- C^2\left[ (f^\prime)^2 (3-4\eta^2_\e v + (\eta^2_\e v))+ 2 f^\prime f (\eta^2_\e v)^2 - 2 f^{\prime\prime}f(1-\eta^2_\e v)\right] v^2 \\[8pt] 
&& \displaystyle - nC(1-\eta^2_\e v)^{3/2}f^\prime f \eta^\prime_\e (\eta_\e v^{1/2}) e^{u(q_0)} H(e^{u(q_0)}q_0) v^{3/2}\\[8pt] 
&&\displaystyle + \left[- (n-1)f^2  (1-\eta^2_\e v)^2 -C^2 (f^\prime)^2 f (\eta^2_\e v) (1-\eta^2_\e v) - nCf f^\prime (1-\eta^2_\e v)^2 (\eta^2_\e v) e^{u(q_0)} H(e^{u(q_0)}q_0) \right.\\[8pt]
&&\displaystyle  \left.\ \ \ \  -n f^2 (1-\eta^2_\e v)^{5/2} e^{u(q_0)}  \frac{\partial}{\partial \lambda} \left(\lambda H(\lambda q)\right){\big|_{\lambda=e^{u(q_0)}}}\right] v\\[8pt]
&&\displaystyle - nf^2 (1-\eta^2_\e v)^{5/2} e^{2u(q_0)}  \nabla_0^T H(e^{u(q_0)}q_0) \cdot \nabla u\\[8pt]
&&\displaystyle + nC f^\prime f (1-\eta^2_\e v)^{2}  \left(1-\sqrt{1-\eta^2_\e v} \ e^{u(q_0)} H(e^{u(q_0)}q_0)\right) \leq 0,
\end{array}
\end{equation}
and the proof of the general inequality is complete. Now we prove \eqref{eqlemmagradient}. Taking $f\equiv1$, and dividing \eqref{eqT25} by $(1-\eta^2_\e v)^{2}$, we get that
\begin{equation}\label{eqlemmagradientprelim}
\left[- (n-1) -n (1-\eta^2_\e v)^{1/2} e^{u(q_0)}  \frac{\partial}{\partial \lambda} \left(\lambda H(\lambda q)\right){\big|_{\lambda=e^{u(q_0)}}}\right] v- n (1-\eta^2_\e v)^{1/2} e^{2u(q_0)}  \nabla_0^T H(e^{u(q_0)}q_0) \cdot \nabla u\ \leq 0.
\end{equation}
Assume that $v(q_0)\neq0$ (otherwise $v\equiv 0$ and there is nothing to prove). To conclude the proof it remains to estimate the term $\nabla_0^T H(e^{u(q_0)}q_0) \cdot \nabla u$. 
To this end, recalling the notations used in the proof of Lemma \ref{lemmaunifelliptic}, we define $\tilde h \in \R^n$ as the vector whose $i$-th component is $\tilde h_i:=\nabla_0^T H(e^{u(q_0)}q_0) \cdot \frac{\partial_i}{\|\partial_i \|_{n+1}}$, $i=1,\ldots,n$, where $\partial_i=\frac{\partial \phi}{\partial y_i}(F(q_0))$.  Then, by construction and Cauchy-Schwartz's inequality we have
\begin{equation}\label{eqstimanormahtild}
\|\tilde h\|^2_n = \sum_{i=1}^n \left(\nabla_0^T H(e^{u(q_0)}q_0) \cdot \frac{\partial_i}{\|\partial_i \|_{n+1}}\right)^2 \leq n \|\nabla_0^T H(e^{u(q_0)}q_0)\|_{n+1}^2. 
\end{equation}
 Now, exploiting \eqref{formulagradrev} we have
$$\nabla_0^T H(e^{u(q_0)}q_0) \cdot \nabla u = \lambda^{-2}\sum_{i=1}^n\frac{\partial \tilde u}{\partial y_i}  \nabla_0^T H(e^{u(q_0)}q_0) \cdot \partial_i = \lambda^{-1}\sum_{i=1}^n\frac{\partial \tilde u}{\partial y_i} \tilde h_i=  \lambda^{-1}(\nabla_0 \tilde u, \tilde h)_n,$$
 and thus, from \eqref{formulatrasfgrad}, \eqref{eqstimanormahtild} we deduce that
\begin{equation}\label{eqgradestfinal}
|\nabla_0^T H(e^{u(q_0)}q_0) \cdot \nabla u|=\lambda^{-1}|(\nabla_0 \tilde u, \tilde h)_n | \leq \lambda^{-1} \| \nabla_0 \tilde u\|_n\|\tilde h\|_n\leq \sqrt{n} |\nabla u|\|\nabla_0^T H(e^{u(q_0)}q_0)\|_{n+1}.
\end{equation}
Finally, combining \eqref{eqlemmagradientprelim}, \eqref{eqgradestfinal} and dividing by $v^{1/2}$, we obtain \eqref{eqlemmagradient}. The proof is complete.

\end{proof}

\begin{remark}\label{remsustimagrad}
Applying the gradient estimate \eqref{eqlemmagradient} to the solutions of $Q_{\e}^t(u)=0$, we obtain
\begin{equation}\label{eq:remgradestrev}
\begin{array}{lll}
&&\displaystyle \left[- (n-1) -nt (1-\eta^2_\e |\nabla u(q_0)|^2)^{1/2} e^{u(q_0)}  \frac{\partial}{\partial \lambda} \left(\lambda H(\lambda q)\right){\big|_{\lambda=e^{u(q_0)}}}\right] |\nabla u(q_0)|\\[12pt]
&&\displaystyle \ \ \ \ \ \ \ \ \ \ \ \ \   - n^{3/2} t (1-\eta^2_\e |\nabla u(q_0)|^2)^{1/2} e^{2u(q_0)} \|\nabla_0^T H(e^{u(q_0)}q_0)\|_{n+1} \ \leq 0.
\end{array}
\end{equation}
Hence, it is not possible, by using only this strategy, to get a uniform bound with respect to $t$ for $|\nabla u|_\infty$ as in \cite{TW}.  In fact here we deal with functions defined on a manifold with negative Ricci curvature, and thus in \eqref{eq:remgradestrev} we have a term $-(n-1)$, while in \cite{TW}, for the sphere, this term has the opposite sign. We also point out that this trouble does not depend on the choice of the auxiliary function in the proof of Proposition \ref{stimagradint}, as enlightened by \eqref{eqT25}, where the leading term $v^3$ has a negative coefficient. 
\end{remark}

\section{Proof of Theorem \ref{mainteo}}

\begin{proof}[Proof of Theorem \ref{mainteo}]
We first observe that by definition $(\Omega,H)$ is admissible with constant $\theta$, and thus in the proof of Theorem \ref{teopseudosol} we can take $\bar\epsilon=\theta$. Therefore, for any $\e \in (0,\theta)$, there exists a solution $u_\e$ of the regularized problem \eqref{eqpseudoHbumps}. 
Let us choose $\e \in (0,\theta)$ sufficiently close to $\theta$ so that
\begin{equation}\label{condiz1}
\frac{\partial}{\partial \lambda} \left(\lambda H(\lambda q) \right)< - \frac{1}{r_1(\e-\e^2/4)^{1/2}}, \  \hbox{for all }\ q \in \overline\Omega, \lambda \in [r_1,r_2],
\end{equation}

\begin{equation}\label{condiz2}
\left\|\nabla_0^T H(X)\right\|_{n+1} < \frac{1-\e}{n^{3/2}r_2^2}, \ X\in \mathcal{C}_{\overline{\Omega}}(r_1,r_2),
\end{equation}
and let $u$ be the solution of the regularized problem \eqref{eqpseudoHbumps}.
 
Let $q_0 \in \overline\Omega$ be the maximum point of $|\nabla u|$, and set $v=|\nabla u(q_0)|^2$. There are only two possibilities: $v < (1-\e)^2$ or $v \geq (1-\e)^2$. In the first case there is nothing to prove, in fact, by definition of $\eta_\e$ we have that $u$ is a solution of Problem \ref{eqHbump} and we are done. Therefore let us assume that $v \geq (1-\e)^2$. We point out that in this case $q_0$ cannot belong to $\partial \Omega$ because by Step 6 of the proof of Theorem \ref{teopseudosol} and since $\e < \theta$, we have that $\sup_{\partial\Omega} |\nabla u|^2 \leq (1 - \theta)^2 < (1-\e)^2$. Hence $q_0 \in \Omega$.
We also observe that $u \in C^{3,\beta}(\overline{\Omega})$, for some $\beta \in (0,\alpha]$. In fact, by Theorem \ref{teopseudosol} we know that $u \in C^{2,\beta}_0(\overline\Omega)$. Thanks to Lemma \ref{lemmaczeroestimates} we know that $\Sigma(u)$ is contained in $\mathcal{C}_{\overline\Omega}(r_1,r_2)$ and being $H \in C^{1,\alpha}(\mathcal{C}_{\overline\Omega}(r_1,r_2))$, $\partial\Omega \in C^{3,\alpha}$, by standard regularity results (see \cite{GT}), we get that $u \in C^{3,\beta}(\overline{\Omega})$. Therefore, we can apply Proposition \ref{stimagradint} and recalling that by definition $1-\eta^2_\e v=1-\eta^2_\e\left(|\nabla u(q_0)\right) |\nabla u(q_0)|^2$, we have
\begin{equation}\label{stimaintgrad}
\begin{array}{lll}
&&\displaystyle\left[- (n-1) -n (1-\eta^2_\e v)^{1/2} e^{u(q_0)}  \frac{\partial}{\partial \lambda} \left(\lambda H(\lambda q)\right){\big|_{\lambda=e^{u(q_0)}}}\right] v^{1/2}\\[12pt]
&&\displaystyle \ \ \ \ \ \ - n^{3/2} (1-\eta^2_\e v)^{1/2} e^{2u(q_0)}  \|\nabla_0^T H(e^{u(q_0)}q_0)\|_{n+1} \ \leq 0,
\end{array}
\end{equation}
but on the other hand, since we are assuming that $v \geq (1-\e)^2$, by definition of $\eta_\e$ we have $1\geq 1-\eta^2_\e v \geq  1 - (1-\e/2)^2$, and in view of \eqref{condiz1}, \eqref{condiz2} we have
\begin{eqnarray*}
&&\left[- (n-1) -n (1-\eta^2_\e v)^{1/2} e^{u(q_0)}  \frac{\partial}{\partial \lambda} \left(\lambda H(\lambda q)\right){\big|_{\lambda=e^{u(q_0)}}}\right] v^{1/2}\\
&& \ \ \ \ \ \ -  n^{3/2} (1-\eta^2_\e v)^{1/2} e^{2u(q_0)}  \|\nabla^T H(e^{u(q_0)}q_0)\|_{n+1}\\[8pt]
&>& \left[ -(n-1) + n \frac{e^{u(q_0)}}{r_1} \frac{(1-\eta^2_\e v)^{1/2}}{(\e - \e^2/4)^{1/2}}\right] v^{1/2} - \frac{e^{2u(q_0)}}{r_2^2}(1-\e)\\
&\geq& \left[ -(n-1) + n \frac{(\e-\e^2/4)^{1/2}}{(\e - \e^2/4)^{1/2}}\right]v^{1/2} -(1-\e)\\
&\geq& (1-\e)- (1-\e)= 0\\
\end{eqnarray*}
and thus we contradict \eqref{stimaintgrad}. Therefore the only possibility is $v<(1-\e)^2$, and by definition of $\eta_\e$ this means that $u$ is a solution of Problem \ref{eqHbump}. Moreover, as proved in Theorem \ref{teopseudosol}, such solution  is the unique solution whose associated radial graph is contained in $\mathcal{C}_{\overline\Omega}(r_1,r_2)$, and this completes the proof.  
\end{proof}

\section{A finer gradient estimate}
In this section we prove an a priori estimate for the gradient of the solutions of
\begin{equation}\label{eqdivformt}
\begin{cases}
\displaystyle -\mathrm{div}_{\H^n}\Big(\frac{\nabla u}{\sqrt{1-|\nabla u|^{2}}}\Big)+\frac{nt}{\sqrt{1-|\nabla u|^{2}}}=nt e^{u}H(e^{u}q)&\text{in }\Omega,\\
|\nabla u|<1&\text{in }\Omega,\\
u=0&\text{on }\partial\Omega,
\end{cases}
\end{equation}
where $t \in [0,1]$. As in Sect. 4 we introduce the function $\nu=\frac{1}{\sqrt{1-|\nabla u|^{2}}}$.

\begin{proposition}\label{fineraprioriestgrad}
Let $\Omega$ be a bounded domain of $\H^n$, let $H \in C^{1}(\mathcal{C}_{\overline\Omega})$, let $r_1,r_2 \in \R$ be such that $r_1\neq r_2$, $0<r_1\leq1 \leq r_2$ and let $\nu_0 >0$ be a positive number. Then, there exists a constant $C=C(r_1,r_2,\nu_0,\Omega,H)>0$ such that for any $t \in [0,1]$, for any solution $u \in C^{3}(\overline{\Omega})$ of \eqref{eqdivformt} satisfying $\log r_1 \leq u \leq \log r_2$ and $\sup_{\partial\Omega} \nu \leq \nu_0$, we have
$$\sup_{\Omega} \nu\leq C.$$
\end{proposition}

\begin{proof}
Let $u \in C^{3}(\overline{\Omega})$ be a solution of \eqref{eqdivformt} satisfying $\log r_1 \leq u \leq \log r_2$ and $\sup_{\partial\Omega} \nu \leq \nu_0$. Clearly $\nu \in C^0(\overline\Omega)$ and we can introduce the differential operator $P_u:C^1(\overline\Omega) \to C^0(\overline\Omega)$ defined by $$P_uw:=\nu\sum_{k=1}^nu_kw_k,$$ where $u_k$, $w_k$ are the covariant derivatives with respect to a orthonormal frame field. Applying $P_u$ to both sides of the equation in \eqref{eqdivformt} and arguing as in Proposition \ref{stimagradint} we deduce that $\nu$ satisfies the following equation:

\begin{equation}\label{eq1fgrad}
\begin{array}{lll}
-\displaystyle \sum_{i,j=1}^n\nabla_i\left(\nu^{-2}f_{ij}\nu_j\right)&+&\nu^{-2}|\nabla \nu|^2 +|\l\nabla u, \nabla \nu\r|^2 +\displaystyle\nu\sum_{i,j,k=1}^n f_{ij} u_{jk}u_{ik}+ \sum_{i,j=1}^n\nu^2 R_{ij} u_i u_j\\ 
&+&  nt \nu \l\nabla u, \nabla \nu\r=\displaystyle   \nu n t \sum_{k=1}^n u_k \nabla_k(e^{u}H(e^{u}q)),
\end{array}
\end{equation}
where $f_{ij}:=\nu \delta_{ij}+ \nu^3 u_iu_j$, $R_{ij}=-(n-1)\delta_{ij}$ is the Ricci curvature tensor of $\H^n$, $i,j=1,\ldots,n$. This relation resemble that appearing in \cite[(4.8)]{Gerhardt}, and it can be proved by direct computation taking into account of the identities $\nu^{-2}=1-|\nabla u|^2$, $\nu_i=\nu^3 \sum_{l=1}^n u_l u_{li}$, $\l\nabla u, \nabla \nu\r=\nu^3\sum_{l,m=1}^nu_lu_m u_{lm}$, $|\nabla \nu|^2=\nu^6\sum_{i=1}^n\left( \sum_{l=1}^n u_l u_{li}\right)^2$, and \cite[(2.6)]{Yau}. In order to estimate the terms appearing in \eqref{eq1fgrad} we first observe that

\begin{equation}\label{eq2fgrad}
\begin{array}{lll}
\displaystyle\sum_{i,j,k=1}^n f_{ij} u_{jk}u_{ik}&=&\displaystyle\sum_{i,j,k=1}^n \left(\nu \delta_{ij}+ \nu^3 u_iu_j\right)u_{jk}u_{ik}\\
&=&\displaystyle \nu \sum_{i,k=1}^n u_{ik}^2 + \nu^3 \sum_{k=1}^n \left(\sum_{i=1}^nu_i u_{ik}\right)^2\\[16pt]
&\geq& \displaystyle \nu |D^2u|^2, 
\end{array}
\end{equation}
where $|D^2u|^2:= \sum_{i,k=1}^n u_{ik}^2$ is the square of the matrix norm of the Hessian. For the term $nt \l\nabla u, \nabla \nu\r$, we write the equation in \eqref{eqdivformt} in non-divergence form as 

\begin{equation}\label{eq5fgrad}
-\nu \Delta u - \l \nabla u, \nabla \nu\r + nt \nu=nt e^uH(e^uq).
\end{equation}

Then, multiplying each side by $nt \nu$, recalling that $\nu\geq1$, $e^uH(e^uq)$ is uniformly bounded with respect to $t$, and using the inequality $|\Delta u| \leq \sqrt{n}|D^2u|$ we deduce that
$$|nt \nu \l\nabla u, \nabla \nu\r | \leq c_1\nu^2 (1+|\Delta u|)\leq c_2\nu^2(1+|D^2u|),$$
for some constants $c_1,c_2>0$ depending on $n$, $r_1$, $r_2$ and $\|H\|_{\infty,\mathcal{C}_{\overline\Omega}}$, but not on $t$. From now on $c_3$, $c_4$, etc. will denote positive constants which do not depend on $t$. Now, if $|D^2u|<c_2 \left(1+\sqrt{1+\frac{1}{c_2}}\right)$ we get that
$|nt \nu \l\nabla u, \nabla \nu\r |\leq c_3\nu^2$, where $c_3$ depends just on $c_2$, and thus $nt\nu  \l\nabla u, \nabla \nu\r \geq  -c_3\nu^2$. On the other hand, if $|D^2u|\geq c_2 \left(1+\sqrt{1+\frac{1}{c_2}}\right)$, by an elementary computation we infer that
$$- c_2\nu^2(1+|D^2u|) + \frac{1}{2}\nu^2 |D^2u|^2\geq 0.$$
Hence, in view \eqref{eq2fgrad} and the previous inequalities we obtain
\begin{equation}\label{eq3fgrad}
\displaystyle\sum_{i,j,k=1}^n f_{ij} u_{jk}u_{ik} +  nt \nu \l\nabla u, \nabla \nu\r \geq - c_4 \nu^2 +  \frac{1}{2}\nu^2 |D^2u|^2.
\end{equation}
Therefore, from \eqref{eq1fgrad}, \eqref{eq3fgrad} we have
\begin{equation}\label{eq4fgrad}
-\sum_{i,j=1}^n\nabla_i\left(\nu^{-2}f_{ij}\nu_j\right)  +|\l\nabla u, \nabla \nu\r|^2 +  \frac{1}{2}\nu^2 |D^2u|^2 \leq c_5 \nu^2 + \nu n t \sum_{k=1}^n u_k \nabla_k(e^{u}H(e^{u}q)).
\end{equation}
Now, writing \eqref{eq5fgrad} as $-\nu \Delta u - \l \nabla u, \nabla \nu\r =nt e^uH(e^uq)-nt\nu$ and squaring, then, by using elementary inequalities we get that
\begin{equation}\label{eq4bfgrad}
\nu^2|\Delta u|^2 - 2\nu |\Delta u| |\l\nabla u, \nabla \nu \r| + |\l\nabla u, \nabla \nu \r|^2\leq 2n^2e^{2u}H^2(e^{u}q)+2 n^2\nu^2.
\end{equation}
Multiplying \eqref{eq5fgrad} by $\nu$, and using $|\Delta u| \leq \sqrt{n}|D^2u|$, we deduce that $\nu | \l\nabla u, \nabla \nu\r | \leq c_5\nu^2(1+|D^2u|)$. Hence, from this, using again $|\Delta u| \leq \sqrt{n}|D^2u|$, and \eqref{eq4bfgrad} we obtain
$$-n\nu^2|D^2u|^2 - 2\sqrt{n}c_5 \nu^2(1+|D^2u|) + |\l\nabla u, \nabla \nu \r|^2\leq 2n^2e^{2u}H^2(e^{u}q)+2 n^2\nu^2,$$
and thus by elementary computations we deduce that
\begin{equation}\label{eq5bfgrad}
-c_6 \nu^2|D^2u|^2 + |\l\nabla u, \nabla \nu \r|^2\leq 2n^2e^{2u}H^2(e^{u}q)+ c_7 \nu^2.
\end{equation}
Therefore, dividing \eqref{eq5bfgrad} by $C:=2c_6+1$ and summing with \eqref{eq4fgrad} we deduce the following
 \begin{equation}\label{eq6fgrad}
 \begin{array}{lll}
&&\displaystyle-\sum_{i,j=1}^n\nabla_i\left(\nu^{-2}f_{ij}\nu_j\right)  +(1+2c_*)|\l\nabla u, \nabla \nu\r|^2 +  c_*\nu^2 |D^2u|^2\\
& \leq&\displaystyle c_8 \nu^2 + c_8 e^{2u}H^2(e^{u}q)+ \nu n t \sum_{k=1}^n u_k \nabla_k(e^{u}H(e^{u}q)),
\end{array}
\end{equation}
where $c_{*}=\frac{1}{2}-\frac{c_6}{2c_6+1}>0$ does not depend on $t$. From \eqref{eq6fgrad}, by arguing as in \cite[Theorem 4.1]{Gerhardt}, we can conclude the proof. In fact, using Stampacchia's truncation method (for the details see the Appendix in \cite{Gerhardt}), multiplying \eqref{eq6fgrad} with
$$\psi_l:=\nu \max\{\nu-l,0\}, \ \ l\geq \nu_0,$$
and integrating by parts we deduce the following 

\begin{equation}\label{eqstimal2n}
 \sup_\Omega  \nu \leq \nu_0 + c_9 (1+|\nu|_{2n,\Omega}^3),
\end{equation}
where $c_9>0$ is a constant depending on $n$, $\Omega$, $r_1$, $r_2$ and $\|H\|_{\infty,\mathcal{C}_{\overline\Omega}}$ but not on $t$, and $|\cdot|_{p,\Omega}$ denotes the standard $L^p$-norm. Therefore, in order to conclude the proof it suffices to prove a uniform estimate for the $L^{2n}$-norm of $\nu$ with respect to the parameter $t$. To this end, recalling that $\nu \geq 1$, and that $e^uH(e^u)$ is uniformly bounded by a constant depending only on $r_1,r_2, \|H\|_{\infty, \mathcal{C}_{\overline	\Omega}}$, we can rewrite the right-hand side of \eqref{eq6fgrad} in a simpler way 

 \begin{equation}\label{eq7fgrad}
-\sum_{i,j=1}^n\nabla_i\left(\nu^{-2}f_{ij}\nu_j\right)  +(1+2c_*)|\l\nabla u, \nabla \nu\r|^2 +  c_*\nu^2 |D^2u|^2\\
 \leq \displaystyle c_{10} \nu^2 + \nu n t \sum_{k=1}^n u_k \nabla_k(e^{u}H(e^{u}q)).
\end{equation}
Now, let $p\geq 2$ any fixed real number, let $\lambda>0$ be a real number to be chosen later and multiply \eqref{eq7fgrad} with
$$\rho_l:=\nu_l^pe^{\lambda u},$$ where $\nu_l := \max\{\nu-l,0\}$ and $l$ is any fixed number such that $l\geq \nu_0$. Since $\nu_l^pe^{\lambda u} \in H_0^{1,q}(\Omega)$, for any $q \in [1,+\infty[$, we can integrate by parts and thus we obtain
 
 \begin{equation}\label{eq8fgrad}
 \begin{array}{lll}
&&\displaystyle p\sum_{i,j=1}^n \int_\Omega \nu^{-2 } f_{ij}\nu_j \nu_i \nu_l^{p-1}e^{\lambda u} + \lambda\sum_{i,j=1}^n \int_\Omega \nu^{-2 } f_{ij}\nu_j u_i \nu_l^{p}e^{\lambda u}\\[12pt]
& +&\displaystyle(1+2c_*)\int_\Omega |\l\nabla u, \nabla \nu\r|^2  \nu_l^{p} e^{\lambda u} +  c_*\int_\Omega\nu^2 |D^2u|^2  \nu_l^{p}e^{\lambda u}\\[6pt]
& \leq&\displaystyle c_{11} \int_\Omega \nu^2  \nu_l^{p}e^{\lambda u}+ c_{11}(p+1) \int_\Omega  \nu  \nu_l^{p-1}  |\l\nabla u, \nabla \nu\r| e^{\lambda u} + c_{11} \lambda  \int_\Omega \nu  \nu_l^{p} |D^2u|e^{\lambda u}.
\end{array}
\end{equation}
Now let us observe that

 \begin{equation}\label{eq9fgrad}
 \begin{array}{lll}
\displaystyle \sum_{i,j=1}^n f_{ij}\nu_i \nu_j &=&\displaystyle   \sum_{i,j=1}^n (\delta_{ij} \nu +\nu^3 u_i u_j)\nu_i \nu_j\\[12pt]
&=&\displaystyle \nu |\nabla \nu|^2 +\nu^3 |\l\nabla u, \nabla \nu\r|^2\\[6pt]
&\geq& \nu^3 |\l\nabla u, \nabla \nu\r|^2.
\end{array}
\end{equation}
In addition, by direct computation we have
 \begin{equation}\label{eq10fgrad}
 \begin{array}{lll}
\displaystyle \lambda\sum_{i,j=1}^n \int_\Omega \nu^{-2 } f_{ij}\nu_j u_i \nu_l^{p}e^{\lambda u}&=&\displaystyle  \lambda\sum_{i,j=1}^n \int_\Omega (1-|\nabla u|^2)\nu (\delta_{ij} \nu +\nu^3 u_i u_j)\nu_j u_i \nu_l^{p}e^{\lambda u}  \\[16pt]
&=&\displaystyle \lambda  \int_\Omega \l \nabla u, \nabla \nu\r \nu \nu_l^p e^{\lambda u}.
\end{array}
\end{equation}
Furthermore, fixing a large constant $C_1$, splitting the domains of the integrals in two parts: $|\l\nabla u ,\nabla \nu\r|\leq C_1$ and $|\l\nabla u ,\nabla \nu\r| > C_1$, then, by elementary computations it follows that for a suitable large constant $c_{12}>0$ it hods that
 
 \begin{equation}\label{eq11fgrad}
 \begin{array}{lll}
&&\displaystyle c_*\int_\Omega  \nu_l^{p}  |\l\nabla u, \nabla \nu\r|^2 e^{\lambda u} - c_{11}(p+1) \int_\Omega  \nu  \nu_l^{p-1}  |\l\nabla u, \nabla \nu\r| e^{\lambda u}\\ 
&\geq& \displaystyle - c_{12} \int_\Omega  \nu_l^{p}  e^{\lambda u} - c_{12} \int_\Omega  \nu  \nu_l^{p-1}  e^{\lambda u}.
\end{array}
\end{equation}
Again by elementary considerations we obtain the further estimate
 \begin{equation}\label{eq12fgrad}
 \begin{array}{lll}
&&\displaystyle c_*\int_\Omega\nu^2 |D^2u|^2  \nu_l^{p}e^{\lambda u}- c_{11} \lambda  \int_\Omega \nu  |D^2u| \nu_l^{p}  e^{\lambda u},\\[6pt]
&\geq&\displaystyle-c_{13}\int_\Omega \nu_l^{p}e^{\lambda u}.
\end{array}
\end{equation}
Indeed, since it is always possible to find a constant $c_{11}>0$ such that $c_* x^2 - c_{11} \lambda x + c_{13}>0$ for all $x\geq0$, then, taking $x=\nu |D^2u|$ we obtain the desired inequality.
Therefore, from \eqref{eq8fgrad}, and using the estimates \eqref{eq9fgrad}--\eqref{eq12fgrad}, we deduce that
 
 \begin{equation}\label{eq12afgrad}
 \begin{array}{lll}
&&\displaystyle(1+p+c_*)\int_\Omega |\l\nabla u, \nabla \nu\r|^2  \nu_l^{p} e^{\lambda u}\\
& \leq&\displaystyle \lambda  \int_\Omega  |\l\nabla u, \nabla \nu\r| \nu  \nu_l^{p} e^{\lambda u} + c_{14}  \int_\Omega \nu^2  \nu_l^{p} e^{\lambda u} +\underbrace{c_{14}  \int_\Omega \nu_l^{p} e^{\lambda u} + c_{14}  \int_\Omega \nu  \nu_l^{p-1} e^{\lambda u}}_{(I)}.
\end{array}
\end{equation}

Observe that $(I)$ contains only powers of the form $\nu^{a}\nu_l^{b}$ with $a,b\geq 0$ such that $a+b\leq p+1$. From now on we will denote by $I_1$, $I_2$, etc., terms which are finite sums of integrals of the form $c\int_\Omega \nu^{a}\nu_l^{b} e^{\lambda u}$, where $a+b\leq p+1$, $a,b\geq 0$ and $c$ is a constant which does not depend on $t$. The strategy to conclude the proof is to obtain an estimate of the kind
\begin{equation}\label{eqaim}
 \int_\Omega \nu^2  \nu_l^{p} e^{\lambda u} \leq I.
\end{equation}
To this aim, from \eqref{eq12afgrad}, dividing each side by $(p+1+c_*)$ and using the elementary inequality $xy\leq \frac{1}{2}x^2+\frac{1}{2} y^2$, we obtain that
 \begin{equation}\label{eq13fgrad}
 \displaystyle\int_\Omega |\l\nabla u, \nabla \nu\r|^2  \nu_l^{p} e^{\lambda u} \leq\displaystyle \frac{\lambda^2}{(1+p+c_*)^2}  \int_\Omega \nu^2  \nu_l^{p} e^{\lambda u}+ \frac{2c_{12}}{1+p+c_*}  \int_\Omega \nu^2  \nu_l^{p} e^{\lambda u}  + I_1.
\end{equation}
Now, multiplying \eqref{eq5fgrad} by $\varphi=\nu \nu_l^p e^{\lambda u}$, integrating by parts, taking into account that $\nabla \varphi=\nabla \nu \nu_l^p e^{\lambda u} + p \nu \nabla \nu \nu_l^{p-1} e^{\lambda u}+\lambda  \nu \nu_l^p \nabla u$, and $p\geq 2$, we get that

 \begin{equation}\label{eq14fgrad}
\displaystyle\lambda \int_\Omega \nu^2  \nu_l^{p} |\nabla u|^2 e^{\lambda u} \leq\displaystyle c_{15}\int_\Omega \nu^2  \nu_l^{p} e^{\lambda u}+ (p+1) \int_\Omega \nu^2  \nu_l^{p-1} |\l\nabla u, \nabla \nu\r| e^{\lambda u}  + I_2.
\end{equation}
Now, choosing $\lambda>c_{15}$, recalling that $\nu^{-2}=1-|\nabla u|^2$, from \eqref{eq14fgrad} we obtain

 \begin{equation}\label{eq15fgrad}
\lambda \int_\Omega \nu^2  \nu_l^{p} e^{\lambda u}\leq\displaystyle \frac{(p+1)^2}{(\lambda-c_{15})^2} \int_\Omega \nu_l^{p} |\l\nabla u, \nabla \nu\r|^2 e^{\lambda u}  + I_3.
\end{equation}
From the combination of \eqref{eq13fgrad} and \eqref{eq15fgrad}, for a large $\lambda$ such that $$\frac{\lambda^2(p+1)^2}{(p+1+c_*)^2(\lambda-c_{15})^2}+\frac{2c_{12}(p+1)^2}{(p+1+c_*)(\lambda-c_{15})^2}<1$$ it follows that

$$ \int_\Omega \nu_l^{p} |\l\nabla u, \nabla \nu\r|^2 e^{\lambda u} \leq I_4,$$
and then, from this and \eqref{eq15fgrad}, we conclude that
$$ \int_\Omega \nu^2  \nu_l^{p} e^{\lambda u}\leq I_5,$$
which gives the desired inequality \eqref{eqaim}. Therefore, from \eqref{eqaim} and the arbitrariness of $p$ we deduce that $|\nu|_{2n,\Omega}$ is uniformly bounded in $t$ and thus from \eqref{eqstimal2n} we deduce the thesis. The proof is complete.
\end{proof}

\section{Proofs of Theorem \ref{mainteoprob} and Theorem \ref{mainteoprob2}}
The proofs of Theorem \ref{mainteoprob} and Theorem \ref{mainteoprob2} are identical except for a small part and thus we give a unified proof in which at some point we distinguish between the two cases.
\begin{proof}
Let $\alpha$, $r_1,r_2$, $\Omega$ and $H$ be as in the statement of the theorem. Recalling the definition of the operators $\mathcal{Q}^t$, $\hat{\mathcal{Q}^t}$ (see \eqref{operatorcomp}, \eqref{operatorcomp2}), then by the same proof of Lemma \ref{lemmaczeroestimates} it follows that, for any $t \in [0,1]$, if $u \in C^2_0(\overline\Omega)$ is a solution of  $\hat{\mathcal{Q}^t}(u)=0$ and satisfies $|\nabla u|_{\infty,\Omega}<1$ then
\begin{equation}\label{eq1mainteocomp}
\log r_1 \leq u(q) \leq \log r_2, \ \ \hbox{for any} \ q \in \overline\Omega.
\end{equation}
Hence, by definition of $\mathcal{Q}^t$, we have also a uniform bound with respect to $t$ on the $L^\infty$ norm of the solutions of $\mathcal{Q}^t(u)=0$. 
In order to get a uniform bound on the gradient we use Proposition \ref{fineraprioriestgrad}. To this end, in the case of Theorem \ref{mainteoprob} since $\Omega$ satisfies a uniform exterior geodesic condition and $H>0$, then, thanks to Proposition \ref{propadmissextbgball} we have that $(\Omega,H)$ is admissible, and by arguing as in Step 6 of the proof of Theorem \ref{teopseudosol} we obtain that there exists $\theta \in (0,1)$ such that for any $t \in [0,1]$, if $u \in C^2_0(\overline\Omega)$ is a solution of  ${\mathcal{Q}^t}(u)=0$ and satisfies $|\nabla u|_{\infty,\Omega}<1$, then
$$|\nabla u(q)|\leq  1-\theta, \ \hbox{for any} \ q \in \partial\Omega.$$
Indeed, if $|\nabla u|_{\infty,\Omega}<1$ and $u \in C^1(\overline\Omega)$, then, by the same proof of Lemma \ref{lemmaunifelliptic} we get that ${\mathcal{Q}^t}$ is uniformly elliptic in $\Omega$ (when passing to hyperbolic stereographic coordinates) and thus, thanks to the hypotheses on $H$, we can apply \cite[Theorem 10.1]{GT} and argue as in Step 6 of the proof of Theorem \ref{teopseudosol}. In the case of Theorem \ref{mainteoprob2} the proof of this fact is identical and we use directly the hypothesis that $(\Omega,H)$ is admissible without invoking Proposition \ref{propadmissextbgball}.

Since $\Omega$ is of class $C^{3,\alpha}$, $H \in C^{1,\alpha}(\mathcal{C}_{\overline\Omega}(r_1,r_2))$ and thanks to \eqref{eq1mainteocomp} then by standard elliptic regularity theory (see \cite{GT}) any solution $u \in C_0^{2,\alpha}(\overline\Omega)$ of ${\mathcal{Q}^t}(u)=0$ such that $|\nabla u|<1$ in $\overline\Omega$ turns out to be of class $C^{3,\alpha}(\overline\Omega)$. Hence, setting $\nu_0:=  \frac{1}{\sqrt{1-\theta^2}}$, by Proposition \ref{fineraprioriestgrad}, it follows that there exists $\theta_* \in (0,1)$, depending only on $n,r_1,r_2,\nu_0,\Omega,H$ but not on $t$, such that for any solution $u \in C^3(\overline\Omega)$ of ${\mathcal{Q}^t}(u)=0$, satisfying $|\nabla u|<1$ in $\Omega$, it holds that
\begin{equation}\label{eq2mainteocomp}
|\nabla u(q)| \leq 1-\theta_*, \ \hbox{for any} \ q \in \overline\Omega.
\end{equation}
Let us fix $\delta>0$ sufficiently small so that $1-\theta_*+\delta<1$ and consider the set 
$$U:=\{w \in C^{1,\alpha}_0(\overline\Omega);\ |\nabla w|_{\infty,\Omega} \leq 1-\theta_*+\delta\}.$$
Clearly $U$ is a convex and closed subset of $C^{1,\alpha}_0(\overline\Omega)$. We define the map $T:[0,1]\times U \to C^{1,\alpha}_0(\overline\Omega)$, $T(t,w):=u$, where $u$ is the unique solution of
 
\begin{equation*}
\begin{cases}
\displaystyle \sum_{i,j=1}^n \left((1-|\nabla w|^2) \delta_{ij} + w_i w_j\right) u_{ij} = n t (1- |\nabla w|^2) \left(1- \sqrt{1- |\nabla w|^2} e^w  \hat H(e^w q)\right) & \hbox{in} \ \Omega,\\[6pt]
u=0 & \hbox{on} \ \partial\Omega.
\end{cases}
\end{equation*}
We observe that $T$ is well defined. Indeed, for a fixed $w \in U$, considering the linear operator $L_w u:= \sum_{i,j=1}^n \left((1- |\nabla w|^2) \delta_{ij} + w_i w_j\right) u_{ij}$, and arguing as in Step 2 of the proof of Theorem \ref{teopseudosol},  we see that $L_{w,\e}:C^{2,\alpha}_0(\overline\Omega) \to C^{0,\alpha}(\overline \Omega)$ is a bijection. Hence $$T(w)=t L_w^{-1}\left( n(1- |\nabla w|^2) \left(1- \sqrt{1- |\nabla w|^2} e^w  \hat H(e^w q)\right)\right)$$
is defined and we are done.

It is easy to verify that $T$ is continuous and, arguing as in the proof of Step 5 of Theorem \ref{teopseudosol}, we have that $T([0,1]\times U)$ is a relatively compact subset of $C^{1,\alpha}_0(\overline\Omega)$. Moreover $0$ lies in the interior of $U$ and $T(0\times\partial U)\subset U$. To conclude the proof it suffices to prove that if $(t,u) \in [0,1]\times U$ satisfies $T(t,u)=u$ then $u \not\in \partial U$. Indeed, if $T(t,u)=u$ then $u\in C^{2,\alpha}_0(\overline\Omega)$ is a solution of $\hat{\mathcal{Q}^t}(u)=0$ and thus from \eqref{eq1mainteocomp} we have $\mathcal{Q}^t(u)=0$. Then, since $u \in U$ we have $|\nabla u|_{\infty,\Omega}\leq 1-\theta^*+\delta<1$ and thus $\mathcal{Q}^t$ is uniformly elliptic. Therefore by elliptic regularity theory $u \in C^{3,\alpha}(\overline\Omega)$ and thanks to \eqref{eq2mainteocomp} it follows that $|\nabla u|_{\infty,\Omega}\leq 1-\theta_*<1-\theta_*+\delta$, thus $u$ cannot belong to $\partial U$ and we are done.

At the end, from Theorem \ref{LSPTheom} we conclude that there exists $\bar u \in U$ which solves $T(1,\bar u)=\bar u$, i.e., $\bar u$ is a solution of \eqref{eqHbump}. For the uniqueness it suffices to argue as in Step 8 of the proof of Theorem \ref{teopseudosol}. The proof is complete

\end{proof}

\end{document}